\documentclass[11pt]{article}
\usepackage[english]{babel}
\usepackage[latin1]{inputenc}
\usepackage[T1]{fontenc}

\usepackage{amsmath, amsfonts, amssymb, amsthm,graphicx, graphics,psfrag, caption}
\newtheorem*{thm1}{Theorem 1}
\newtheorem*{thm1prime}{Theorem 1'}
\newtheorem*{thm2}{Theorem 2}
\newtheorem{theo}{Theorem}[section]
\newtheorem{lemma}[theo]{Lemma}
\newtheorem{prop}[theo]{Proposition}
\newtheorem{coro}[theo]{Corollary}

\newtheorem{rema}[theo]{Remark}
\newtheorem{claim}[theo]{Claim}

\newcommand{\R}{\mathbb{R}}
\newcommand{\N}{\mathbb{N}}
\newcommand{\T}{\mathbb{T}}

\renewcommand{\P}{\mathbb{P}}
\newcommand{\E}{\mathbb{E}}

\newcommand{\rmP}{\mathrm{P}}
\newcommand{\rmE}{\mathrm{E}}
\newcommand{\rmPo}{\mathrm{P}^{\omega}}

\newcommand{\rmPot}{\mathrm{P}^{\omega,t}}

\newcommand{\Q}{\mathbb{Q}}

\newcommand{\M}{\mathcal{M}}
\newcommand{\cN}{\mathcal{N}}

\newcommand{\cA}{\mathcal{A}}
\newcommand{\fm}{\mathfrak{M}}
\newcommand{\cK}{\mathcal{K}}

\newcommand{\cB}{\mathcal{B}}

\newcommand{\cP}{\mathcal{P}}
\newcommand{\mi}{M^\infty}

\newcommand{\bB}{\mathbb{B}}

\newcommand{\ind}{\mathbf{1}}
\newcommand{\la}{\langle}
\newcommand{\ra}{\rangle}
\newcommand{\e}{\varepsilon}

\title{The spatial Lambda-Fleming-Viot process: \\ an event-based construction and a lookdown representation}
\author{Amandine V\'eber and Anton Wakolbinger}
\date{\small \today}

\begin{document}
\title{The spatial Lambda-Fleming-Viot process: \\ an event-based construction and a lookdown representation}
\author{A. V\'eber\thanks{AV supported in part by the \emph{chaire Mod\'elisation Math\'ematique et Biodiversit\'e} of Veolia Environnement-\'Ecole Polytechnique-Museum National d'Histoire Naturelle-Fondation X and by the ANR project MANEGE (ANR-09-BLAN-0215).}
\\CMAP - \'Ecole Polytechnique\\Route de Saclay\\91128 Palaiseau Cedex\\ France\\~\\ A. Wakolbinger\thanks{AW acknowledges support of the DFG Priority Programme \emph{Probabilistic Structures in Evolution}, SPP 1590}\\Institut f\"ur Mathematik, Goethe-Universit\"at\\ Robert-Mayer-Str. 10\\60325 Frankfurt am Main
\\Germany} \maketitle

\begin{abstract}
We construct a measure-valued equivalent to the spatial $\Lambda$-Fleming-Viot process (SLFV) introduced in \cite{ETH2008}. In contrast with the construction carried out in \cite{ETH2008}, we fix the realization of the sequence of reproduction events and obtain a \emph{quenched} evolution of the local genetic diversities. To this end, we use a particle representation which highlights the role of the genealogies in the attribution of types (or alleles) to the individuals of the population. This construction also enables us to clarify the state-space of the SLFV and to derive several path properties of the measure-valued process as well as of the labeled trees describing the genealogical relations between a sample of individuals. We complement it with a look-down construction which provides a particle system whose empirical distribution at time $t$, seen as a process in $t$, has the law of the \emph{quenched} SLFV. In all these results, the facts that we work with a fixed configuration of events and that reproduction occurs only locally in space introduce serious technical issues that are overcome by controlling the number of events occurring and of particles present in a given area over macroscopic time intervals.


\end{abstract}


\section{Introduction and motivation}
The spatial $\Lambda$-Fleming-Viot process was introduced in \cite{ETH2008} and \cite{BEK2010} as a new framework for studying the genetic evolution of a population distributed over a continuous space.
In these papers, two-dimensional spaces such as $\mathbb R^2$ or the 2-dimensional torus were chosen since this seems natural for modeling the habitats of biological populations. The model, however, allows to consider any domain $\T\subset \R^d$ as a geographical space, and any compact metric space $\mathbb K$ of types (or `alleles'). We refer to \cite{BEH2009}, \cite{BEV2010},  \cite{BEV2011}, \cite{EV2012} and \cite{BEV2012} for some of the first studies of this process and its relations to other models.

One noteworthy difference with previous models of structured populations is the following:
while in the latter the reproduction mechanism can be phrased in terms of clocks carried by the individuals and deciding of reproduction and death times, in the spatial $\Lambda$-Fleming-Viot process it is driven by a Poisson point process of extinction/recolonization events (of various {\em ranges} $r > 0$ and \emph{impacts} $u \in (0,1]$). For brevity, we shall address such a configuration of events as an {\em environment}. Once an environment $\omega = \{(t_i,z_i,r_i,u_i),\, i\in I\}$ has been fixed, local populations change only at the times $t_i$ of an event, and only within the corresponding ball $\mathbb B(z_i,r_i)\subset \T$. During such an event, a parent is chosen uniformly at random within $\mathbb B(z_i,r_i)$ and a fraction $u_i$ of the local population there is replaced by offspring of that individual in such a way that the total mass of the (continuum) population in this region remains constant. These offspring inherit the allele of their parent, which implies a jump in the local allele frequency at that time. As we shall see, the Poisson point process formulation also enables us to model the ancestral process of a sample of individuals, which, thanks to an inherent assumption of neutrality, is a process \emph{independent of types}: backwards-in-time, all ancestral lineages which are in the geographical area of an event are chosen to be (or not to be) affected by the event in an i.i.d. manner, that is by tossing a coin with success probability $u$.
The affected lineages then coalesce into a single ancestor whose position is uniformly distributed over the area of the event, and the ancestral process starts again from this new configuration. We shall assume that the random environment is a Poisson process with a possibly infinite intensity.

In a non-spatial situation (i.e. with $\mathbb T$ consisting of a single element only), the random environment reduces to a Poisson process of so-called $u$-{\em mergers},
with intensity measure $dt \otimes \nu(du) = dt \otimes u^{-2} \Lambda(du)$.
In  \cite{BLG2003}, the resulting forward evolution of  type frequencies was analysed as a measure-valued process and called $\Lambda$-Fleming-Viot process.

The spatial $\Lambda$-Fleming-Viot process (SLFV, in short) was formally constructed in \cite{BEV2010} using a method of Evans \cite{EVA1997}. This more analytic approach, which we re-describe  in Section~\ref{subs:previous def}, relies on a characterization of the semigroup of the SLFV through a family of duality relations with the genealogical process mentioned above.

Our first motivation for this work is to relate the \emph{forwards} and \emph{backwards} evolutions in a more detailed manner, by first defining a reproduction model conditioned on the environment and then by attributing the types thanks to an explicit use of the genealogical relations between the individuals of the population. This construction, which leads to Theorem~1, provides a more probabilistic approach to the definition of the process and enables us to include a mutation mechanism (which was not incorporated in  the analytic construction and in fact seems less tractable there). In order to further disentangle the different sources of randomness acting on our population, in Section~\ref{genealogy} we extend the environment $\omega$ by adding a fifth coordinate recording the spatial location of the parent during each event. This allows to define a {\em parental skeleton}, out of which emerge all the genealogical trees relating a sample of individuals. In Section~\ref{flow}, we randomly label the individuals sitting at each site by  numbers from $[0,1]$, and extend the environment by a sixth coordinate that records the label of the parent during each reproduction event (in the spirit of Bertoin and Le Gall's flows of bridges, see \cite{BLG2003}). We show that in the absence of mutation, conditionally on this extended environment the SLFV can be seen as a deterministic flow on the field of local allele frequencies which at any time conserves the geographical density of the population (here assumed to be constant over space), and transports types through an accumulation of local reproduction events.  Mutations may distort this flow by acting independently on each ancestral line between the reproduction events.

Our second motivation consists in clarifying the state-space of the SLFV and its topological properties. This is done in Lemma~\ref{lem:topo}. Under the basic integrability assumption of \cite{BEV2010} (see condition \eqref{cond def} below) we obtain some general properties about the process: Lemma~\ref{lem: fv} shows that the \emph{quenched} (i.e., for a fixed environment) SLFV without mutation has c\`adl\`ag paths of finite variation with probability $1$; Proposition~\ref{prop:CDI} tells us that, when the fraction of individuals replaced during an event is always less than $1$, the genealogical process of a countable sample of individuals never \emph{comes down from infinity} (see Section~\ref{section:CDI} for a definition of this property).

Thirdly, we propose a representation of the SLFV as the empirical distribution of the paths of an infinite collection of individuals. This is done in Section~\ref{section:look-down} through a look-down construction, which enables us to approximate the \emph{quenched} spatial $\Lambda$-Fleming-Viot process by the empirical distribution of a large but finite set of individuals, uniformly over compact time intervals. The proof of our Theorem~2 is inspired by that of Theorem~1.1 in \cite{BB+2009}, which considers non-spatial $\Lambda$-Fleming-Viot processes. With the additional spatial structure, the fact that reproduction events happen only locally brings in new technical subtleties.

Since their introduction by Donnelly and Kurtz in \cite{DK1996} and \cite{DK1999}, look-down constructions have proved useful in different contexts. For example, in \cite{BB+2005} it is used to show that the genealogy of an $\alpha$-stable branching process (with $\alpha\in (1,2)$) can be described as a randomly time-changed $\beta$-coalescent. In \cite{PW2006}, it is at the core of the study of the evolving coalescent encoding the genealogy of the whole population in a Wright-Fisher model. Finally, in \cite{LZ2012} it is used to show that a given class of $\Lambda$-Fleming-Viot processes (`coming down from infinity' sufficiently fast) with Brownian mutation have a compact support at any positive time. Such a construction was established for the non-spatial $\Lambda$-Fleming-Viot processes in \cite{BB+2009}. In \cite{DE+2000}, a look-down construction for continuum-sites stepping-stone models was carried out, in which the genealogy consists in Feller processes which coalesce instantly upon meeting. Because the corresponding genealogies are fairly different, so are the details of their proof. Furthermore, they only show the convergence of the finite-dimensional distributions of their look-down process, whereas we prove its convergence in the space of c\`adl\`ag paths. To our knowledge, there are no other constructions in the flavour of \cite{DK1996} in which the spatial structure of the population influences their reproduction. Note however that Etheridge and Kurtz \cite{EK2012} propose a look-down construction of the SLFV in the spirit of \cite{KR2011}, and use it in particular to show that the finite-density population model introduced in \cite{BEH2009} converges towards the \emph{annealed} SLFV as the density of individuals tends to infinity.

In a more applied point of view, the different constructions carried out in Sections~\ref{section:exist} and \ref{section:look-down} allow us to interpret a `real' population with reasonably large density as a  Poissonization of the measures describing the state of the theoretical `infinite' population. That is, at any time $t$ we can think of our `real' population as being the result of the evolution of the individuals with level less than some quantity $N$ in the look-down construction (with $N$ not too small) between times $0$ and $t$, which form a locally finite Poissonian cloud of trajectories in $\T\times \mathbb K$. Furthermore, the definition of the SLFV in a \emph{fixed} environment may serve in reconstructing the main features of the history of a population. Indeed, the genetic diversity observed at a given locus within a sample of individuals gives us some information about quantities such as the local rate of coalescence, which depend essentially on the \emph{annealed} evolution. We refer to \cite{BEKV2012} for more on inference issues related to the SLFV. However, if we now consider a second locus in the genome, or even a fraction of this genome (data which are now accessible through modern sequencing tools), the genetic diversities observed at all these sites are correlated first and foremost by the fact that they all flow through the same environment $\omega$. Hence, any information on $\omega$ brought by the analysis of the diversity at one locus yields some constraints on the genetic diversity at other loci. A good understanding of the \emph{quenched} SLFV is thus necessary to get a handle on the correlations between different parts of the genome and to devise statistical methods based on the available data to detect the presence of atypically big events acting as
major catastrophes,
or of natural selection acting on a given gene, or of any other kind of deviations from the hypothesis of neutral evolution made in this work.

The paper is laid out as follows. In Section~\ref{section:introductory}, we describe the ``kernel-valued'' version of the spatial $\Lambda$-Fleming-Viot process introduced in \cite{ETH2008} and state our main result, namely the existence and uniqueness of the \emph{quenched} SLFV with mutation taking values in the space of Radon measures on $\T\times \mathbb K$. We also present the two extensions of the environment discussed above. Theorem~1 is proved in Section~\ref{section:exist} by constructing the two-parameter semigroup of transition probabilities conditionally on the configuration $\omega$ of events, for almost every $\omega$.
While this construction focusses on fixed times $s<t$, in Section~\ref{section:look-down} we materialize the dynamics of the evolution by a look-down construction directly defined for all times. Though these two steps could in fact be summed up into a single one, we chose to keep them separated to put forward the role of the Poisson structure and of the exchangeability of the population in the construction of the SLFV, and because the addition of levels is somehow an elaboration on the ideas used in the first construction. In Section~\ref{section:CDI}, we use the same kind of arguments as in the previous sections to address the question of \emph{coming down from infinity} for the genealogical process.


\section{The  spatial $\Lambda$-Fleming-Viot process and its sources of randomness}\label{section:introductory}
\subsection{The kernel-valued spatial $\Lambda$-Fleming-Viot process}\label{subs:previous def}
Let us start by describing the first construction of the SLFV, carried out in \cite{ETH2008} and \cite{BEV2010}.  As the geographical space in which the population evolves we fix a domain $\T\subset \R^d$, and we write $\mathbb K$ for the compact space of all possible types. For a locally compact metric space (with $E=\T\times \mathbb K$ or $E=\mathbb K$ as the generic examples), let $\M(E)$ (resp., $\M_1(E)$) denote the set of all nonnegative Radon (resp., probability) measures on $E$, and let $\Xi$ be the quotient of the space of all Lebesgue-measurable maps (or {\em probability kernels}) $\rho:\T\rightarrow \M_1(\mathbb K)$ by the equivalence relation
$$
\rho\sim \rho' \qquad \Leftrightarrow \qquad \hbox{Vol}\big(\big\{x\in \T\,:\, \rho(x, .)\neq \rho'(x, .)\big\}\big)=0.
$$
In other words, two maps $\rho$ and $\rho'$ are in the same equivalence class iff the measures $dx\,\rho(x,d\kappa)$ and $dx\,\rho'(x,d\kappa)$ on $\T\times \mathbb K$ are equal. We shall give a mathematically equivalent description of $\Xi$ at the beginning of the next subsection.

To specify the dynamics, let $\mu$ be a $\sigma$-finite measure on $(0,\infty)$ and $\{\nu_r,\, r>0\}$ be a collection of probability measures on $[0,1]$ such that the map $r\mapsto \nu_r$ is measurable with respect to $\mu$. Let $\Pi$ be a Poisson point process on $\R\times \T \times (0,\infty)\times [0,1]$ with intensity measure $dt\otimes dz\otimes \mu(dr)\nu_r(du)$. The random point configuration $\Pi$ acts as a \emph{random environment} in which the population evolves: if $(t,z,r,u)\in \Pi$, then at time $t$ a {\em reproduction event} occurs within the closed ball $\mathbb B(z,r)\subset \T$, that is
\begin{itemize}
\item a location $y$ is sampled uniformly at random in $\mathbb B(z,r)$ and a type $k$ is chosen according to $\rho_{t-}(y,d\kappa)$ (equivalently, a type is sampled according to the mean type distribution in $\mathbb B(z,r)$ just before the event),
\item for every $x\in \mathbb B(z,r)$, $\rho_t(x,d\kappa):= (1-u)\rho_{t-}(x,d\kappa) + u\delta_k$.
\end{itemize}
The value of $\rho$ outside the ball remains unchanged, that is $\rho_t(x,d\kappa)=\rho_{t-}(x,d\kappa)$ for every $x\notin \bB(z,r)$.

The existence and uniqueness of a Markov process evolving according to this dynamics is proven in Section~4 of \cite{BEV2010} under the condition
\begin{equation}\label{cond def}
\int_0^{\infty}\int_0^1 ur^d\nu_r(du)\mu(dr)<\infty.
\end{equation}
(In fact this is shown there only for $\T=\R^2$, but the proof is identical in all dimensions and for any domain in $\R^d$.) The process is characterized by a family of duality relations, based on the following very simple idea. In the absence of mutation, the type of the parent is  transmitted to all its offspring. Thus, to know the types of a few individuals sampled at time $t$, it suffices to go back to their ancestors at some reference time (say, $0$) in the past, and to check which types they carried. Now, some of the individuals alive at time $t$ may share a common ancestor at time $0$, which introduces some correlations between their types. As a consequence, the natural dual process to consider is that tracing the ancestral relations between a sample of individuals, from time $t$ on and back into the past.

Let us thus imagine what the ancestral process of a sample of individuals should look like under the prescribed evolution. Here we forget about types, and for any $k\in \N$, we represent the genealogical relations of a sample of $k$ individuals (labelled by $1,\ldots,k$) by a process $(\cA_h)_{h\geq 0}$ with values in the set of all partitions of $\{1,\ldots,k\}$ whose blocks are marked by an element of $\T$. In words, $h =0$ corresponds to the time at which our individuals are sampled and for any $h \geq 0$, each block of $\cA_h$ contains the labels of all of these individuals who share a common ancestor $h$ units of time back in the past. The mark of a block records the spatial location of the corresponding ancestor at that time. When convenient, we shall write
\begin{equation}\label{notation A}
\cA_h = \big\{(B_h^1,\xi_h^1),\,\ldots\,,(B_h^{N_h},\xi_h^{N_h})\big\},
\end{equation}
where $\xi^i_h$ is the mark of the block $B^i_h$ and $N_h$ is the number of distinct ancestors at time $h$ in the past. See Figure \ref{ancestries}  (with $h = t-s$), where $k=4$, $N_h = 2$, $B_h^1 = \{1,2,3\}, B_h^2 = \{4\}$.

As a start, suppose that we sample a single individual and trace back where its ancestors were. Since forwards-in-time the population evolves only at the epochs of events of $\Pi$, so does the spatial location of an ancestral lineage. Going back by a sufficient amount of time, one will encounter the first event $(t',z',r',u')$ in the past in which our individual was not only in the area $\mathbb B(z',r')$ hit at that time, but was also part of the fraction $u'$ of the local population replaced by offspring of the elected parent. (Below we shall argue that condition \eqref{cond def} guarantees that the time it takes back to this first event in the past is a.s. strictly positive.)  Hence, at the time $h= t-t'$ of this event, the ancestral line of our individual jumps to the location of its parent, which is by construction uniformly distributed over $\mathbb B(z',r')$. Using the same reasoning, we can then find an earlier event $(t'',z'',r'',u'')$ during which this parent was born. The ancestral lineage at that time jumps to a location uniformly distributed over $\mathbb B(z'',r'')$ and stays there until the time of the event in its past during which it was born, and so on. Now, observe that the sequence of events  experienced by the ancestral lineage backwards-in-time is governed by the law of $\Pi$. Consequently,  the rate at which the lineage jumps can be computed directly in terms of $\mu$ and $\nu_r$: once the radius $r$ and the impact $u$ of the event have been chosen, the volume of centres such that the current location $\xi^i_h$ of a lineage belongs to the range of the event is $\mathrm{Vol}(\mathbb B(\xi^i_h,r)\cap \T)$ and the probability that the lineage belongs to the set of newborns is $u$. Hence, the time to wait before the next jump is exponentially distributed with parameter
$$
\int_0^{\infty}\int_0^1 u\,\mathrm{Vol}\big(\mathbb B(\xi^i_h,r)\cap \T\big)\,\nu_r(du)\mu(dr),
$$
and at that time $\xi^i$ jumps to a new location uniformly distributed over the area of the corresponding event. Under the condition stated in (\ref{cond def}), this jump rate is bounded by a constant independent of the current location of the lineage.

Now think of several individuals sampled at distinct locations. When a past event comes up  that covers at least one of them, each of the lineages within its range has a probability $u$ of being an offspring of the parent chosen, independently of the others. Then, all the lineages which were born in this event trace back to the same parent and therefore merge at that time into a single ancestor (i.e., the corresponding blocks of $\cA_h$ merge into a single block) whose location is uniformly distributed over the area of the event. The other lineages, outside the ball or inside but not within the pool of offspring, remain unaffected. Because nothing happens between the events in $\Pi$ that affect the genealogy, the description of $\cA$ as a Markovian system of coalescing blocks whose marks evolve like a family of correlated jump processes is complete.

There remains to find an appropriate set of test functions. To this end, let $C(E)$ (resp., $C_c(E)$) stand for the set of all continuous (resp. continuous and compactly supported) functions on the space $E$. Let also $\la m,f\ra$ denote the integral of the function $f$ against the measure $m$. If $k\in \N$, $F\in C_c(\T^n)$ and $g_1,\ldots,g_k\in C(\mathbb K)$, let us define the function $G_{\bf g}$ by
\begin{equation}\label{def Gg}
G_{\bf g}(\kappa_1,\ldots,\kappa_k) := \prod_{\jmath=1}^k g_\jmath(\kappa_\jmath)
\end{equation}
and the function $I_k(\cdot\, ; F,g_1,\ldots,g_k)\in C(\Xi)$ by
\begin{equation}\label{def In}
I_k(\rho\, ; F,g_1,\ldots,g_k) := \int_{\T^k}F(x_1,\ldots,x_k)\,\Big\la \bigotimes_{1\leq \jmath\leq k}\rho(x_\jmath,d\kappa_\jmath),\, G_{\bf g} \Big\ra\, dx_1\ldots dx_k.
\end{equation}
A slight modification of Lemma 4.1 in \cite{BEV2010} shows that the family of functions of the form $I_k(\cdot\, ;F,g_1,\ldots,g_k)$, with $k$, $F$ and $g_1,\ldots,g_k$ as above, is dense in $C(\Xi)$. As a consequence, they constitute a relevant set of test functions.

We can now state the duality relations upon which the characterization of the kernel-valued SLFV depends. In Theorem~4.2 of \cite{BEV2010}, the SLFV is defined as the unique $\Xi$-valued Hunt process $(\rho_t)_{t\geq 0}$ satisfying: for every $\rho_0\in \Xi$, $t\geq 0$, $k\in \N$, $F\in C_c(\T^k)$ and $g_1,\ldots,g_k\in C(\mathbb K)$,
\begin{equation}\label{previous dual}
\E_{\rho_0}\big[I_k(\rho_t\,; F,g_1,\ldots,g_k)\big] = \int_{\T^k} F(x_1,\ldots,x_k)\, \rmE_{\wp_k({\bf x})}\!\bigg[\Big\la \bigotimes_{1\leq \iota \leq N_t}\rho_0(\xi^\iota_t,d\kappa_\iota),\prod_{\iota=1}^{N_t}\prod_{\jmath\in B_t^\iota}g_\jmath(\kappa_\iota)\Big\ra\bigg]dx_1\ldots dx_k,
\end{equation}
where $\E$ denotes the expectation w.r. to the process $(\rho_t)_{t\geq 0}$, $\rmE$ is the expectation w.r. to the process $\cA$, randomized over the environment (we shall be more precise on this point in the next section), and
\begin{equation}\label{def pn}
\wp_k({\bf x}) := \big\{(\{1\},x_1),\,\ldots\,,(\{k\},x_k)\big\}.
\end{equation}
The idea behind (\ref{previous dual}) is precisely that expounded above: the type of an individual living at time $t$ is given by the type distribution $\rho_0(\xi^1_t,\cdot)$ at the location of its ancestor at time $0$, and if several individuals have a common ancestor at time $0$ (i.e., $\mathrm{Card}(B_t^\iota)>1$ for some $\iota$), they all share the same type $\kappa_\iota$ with law $\rho_0(\xi_t^\iota,\cdot)$.

\subsection{The SLFV as a measure-valued process in a random environment}\label{subs:new def}
Let $\M_{\lambda}$ be the space of all nonnegative Radon measures on $\T\times \mathbb K$ whose `spatial' marginal is equal to Lebesgue measure $\lambda$ on $\T$.
Using a well-known disintegration theorem (see e.g. \cite{KAL2002}, p.561), it is not difficult to show that $\Xi$ is in one-to-one correspondence with $\M_{\lambda}$, this correspondence being
\begin{align} \label{corr}m(dx,d\kappa) = dx\, \rho(x,d\kappa).
\end{align}
(Here and below, we sometimes shortly write $dx$ instead of $\lambda(dx)$.) We shall endow $\M_{\lambda}$ with the topology $\mathcal{T}_v$ of vague convergence and the associated Borel $\sigma$-field, and shall use $\M_{\lambda}$ as the state space for the measure-valued SLFV. Using the correspondence \eqref{corr}, the function $I_k$ defined in \eqref{def In} reads as
$$
I_k(m\,;F,g_1,\ldots,g_k) := \big\la m^{\otimes k},F\otimes G_{\bf g} \big\ra,
$$
where $G_{\bf g}$ was defined in (\ref{def Gg}). The following lemma guarantees that the set of functions of the form $I_k$ constitutes a wide enough family of tests functions, and reveals in particular that  the convergence of a sequence $(\rho_n)_{n\geq 1}$ of elements of $\Xi$ (as defined in \cite{BEV2010}) is equivalent to the vague convergence of the measures $m_n(dx,d\kappa) := dx\,  \rho_n(x,d\kappa)$ as $n \to \infty$.
\begin{lemma}\label{lem:topo}
\begin{enumerate}
\item[a)] The space $(\M_{\lambda},\mathcal{T}_v)$ is compact.
\item[b)] For  $k \in \mathbb N$, $F \in C_c(\T^k)$, $g_1, \ldots, g_k \in C(\mathbb K)$, the function $m \mapsto I_k(m\,;F,g_1,\ldots,g_k)$ is $\mathcal T_v$-continuous on $\mathcal M_\lambda$.
\item[c)] The linear span of the set of constant functions and of functions of the form $I_k(\cdot\, ;\, F;g_1,\ldots,g_k)$, $k\in \N$, $F\in C_c(\T^k)$ and $g_1,\ldots,g_k\in C(\mathbb K)$, is dense in $C(\M_{\lambda})$.
\end{enumerate}
\end{lemma}

\noindent {\bf Proof of Lemma~\ref{lem:topo}.} \emph{a)} Since the  $\mathcal T_v$-closedness is immediate from the definition of $\M_\lambda$, it suffices to show that  any given sequence $(m_n)$ in $\M_\lambda$ has a subsequence that converges in the vague topology $\mathcal T_v$. To see this, note that $\T\times \mathbb K$ is locally compact and separable, hence by Theorem A2.3 in \cite{KAL2002}, $\M(\T \times
\mathbb K)$ endowed with the vague topology is Polish (and in particular
complete). Therefore, again by the quoted theorem,  $(\M_\lambda, \mathcal T_v)$ is compact.
\\
\emph{b)} This follows from the fact that for any locally compact separable space  $E$ and $k \in \mathbb N$,  the mapping $m \to m^{\otimes k}$ is continuous w.r.to the vague topologies on $\mathcal M(E)$ and $\mathcal M(E^k)$.
\\
\emph{c)} This space is an algebra that separate points in $\M_{\lambda}$ and contains the constants. Since $\M_{\lambda}$ is compact by $a)$, the Stone-Weierstrass theorem gives us the result. \hfill $\Box$

\medskip
In Section~\ref{section:look-down}, we shall need to put a metric on the space $D_{\M_{\lambda}}[0,\infty)$ of all c\`adl\`ag paths with values in $\M_{\lambda}$. We thus recall the following result (see, e.g., Section~1 in \cite{DK1996}).
\begin{lemma}\label{lem: metric}
There exists a sequence $(f_n)_{n\geq 1}$ of uniformly bounded functions in $C_c(\T \times \mathbb K)$ which separates points in $\M_{\lambda}$. Furthermore, if $(g_n)_{n\geq 1}$ is such a sequence, then
$$
d(m,m'):= \sum_{n=1}^\infty \frac{1}{2^n}\, |\la m,g_n\ra - \la m',g_n\ra|,\qquad m,m'\in \M_\lambda
$$
defines a metric for the topology of vague convergence on $\M_\lambda$, while
$$
\Delta((m_t),(m'_t)):= \int_0^\infty e^{-t}d(m_t,m'_t)\, dt
$$
is a metric for the topology of locally uniform convergence on $D_{\M_\lambda}[0,\infty)$.
\end{lemma}

Let us fix again a measure $\mu$ on $(0,\infty)$ and a collection $\{\nu_r,\, r>0\}$ of probability measures on $[0,1]$ satisfying (\ref{cond def}). Let $\Pi$ be a Poisson point process with intensity measure  $dt\otimes dz\otimes \mu(dr)\nu_r(du)$. We write $\P$ for the distribution of $\Pi$, and note that $\P$ assigns full measure to the set $\Omega$ of point configurations $\omega = (t_i,z_i,r_i,u_i)_{i\in I}$ on $\R \times \T\times (0,\infty)\times [0,1]$ with the property that $t_i \neq t_{i'}$ for $i\neq i'$ and that for all $s < t \in \mathbb R$ and each bounded subset $B$ of $\mathbb T$,
$$ \sum_{i: s\le t_i \le t, \, z_i \in B }r_i^d u_i < \infty.$$
Recall the description of the ancestral process $\cA$ given in Section~\ref{subs:previous def}. For $\P$-a.a. $\omega\in \Omega$ and $t\in \R$, let $\rmPot_{\wp_k({\bf x})}$ stand for the law of the system of coalescing jump processes describing the genealogy of $k$ of `individuals' sampled at time $t$ at the locations $\mathbf x = (x_1,\ldots, x_k)$, conditionally on the environment $\omega$. That is, under $\rmPot_{\wp_k({\bf x})}$:
\\\\
\begin{minipage}{9cm}
\begin{itemize}
\item The evolution of $\cA$ starts at time $h=0$ in $\wp_k({\bf x})$ and uses only the events $(t_i,z_i,r_i,u_i)\in \omega$ such that $t_i\leq t$.
\item Whenever one or more lineages belong to the range of an event, each of the lineages within $\bB(z_i,r_i)$ takes part in this event with probability $u_i$ or remains unaffected with probability $1-u_i$, independently of each other. All those lineages which are affected merge into a single lineage whose location is uniformly distributed over $\bB(z_i,r_i)$. Then $\cA$ remains constant equal to its new value $\cA_{t-t_i}$ until the next event of $\omega$ in the past which hits at least one of the lineages {\em and} for which at least one of these lineages takes part in the merging.
\end{itemize}
\end{minipage}
\qquad
\begin{minipage}{7cm}
\psfrag{t}{}
\psfrag{a}{$t$}
\psfrag{m}{$x_1$}
\psfrag{l}{$x_3$}
\psfrag{o}{$x_2$}
\psfrag{n}{$x_4$}
\psfrag{t11}{}
\psfrag{t21}{}
\psfrag{t12}{}
\psfrag{z2}{}
\psfrag{z1}{}
\psfrag{z3}{}
\psfrag{r}{$h$}
\psfrag{s}{$s$}
\psfrag{z}{$\xi_{t-s}^2$}
\psfrag{b}{$\xi_{t-s}^1$}
\psfrag{c}{$\xi_{t-s}^3$}
\psfrag{s1}{}
\psfrag{s2}{$t$}
\psfrag{s}{$s$}
\psfrag{y1}{$x_3$}
\psfrag{y2}{$x_1$}
\psfrag{y3}{$x_4$}
\psfrag{y4}{$x_2$}
\psfrag{x1}{$\xi_{t-s}^1$}
\psfrag{x2}{$\xi_{t-s}^2$}
\includegraphics[width=7cm]{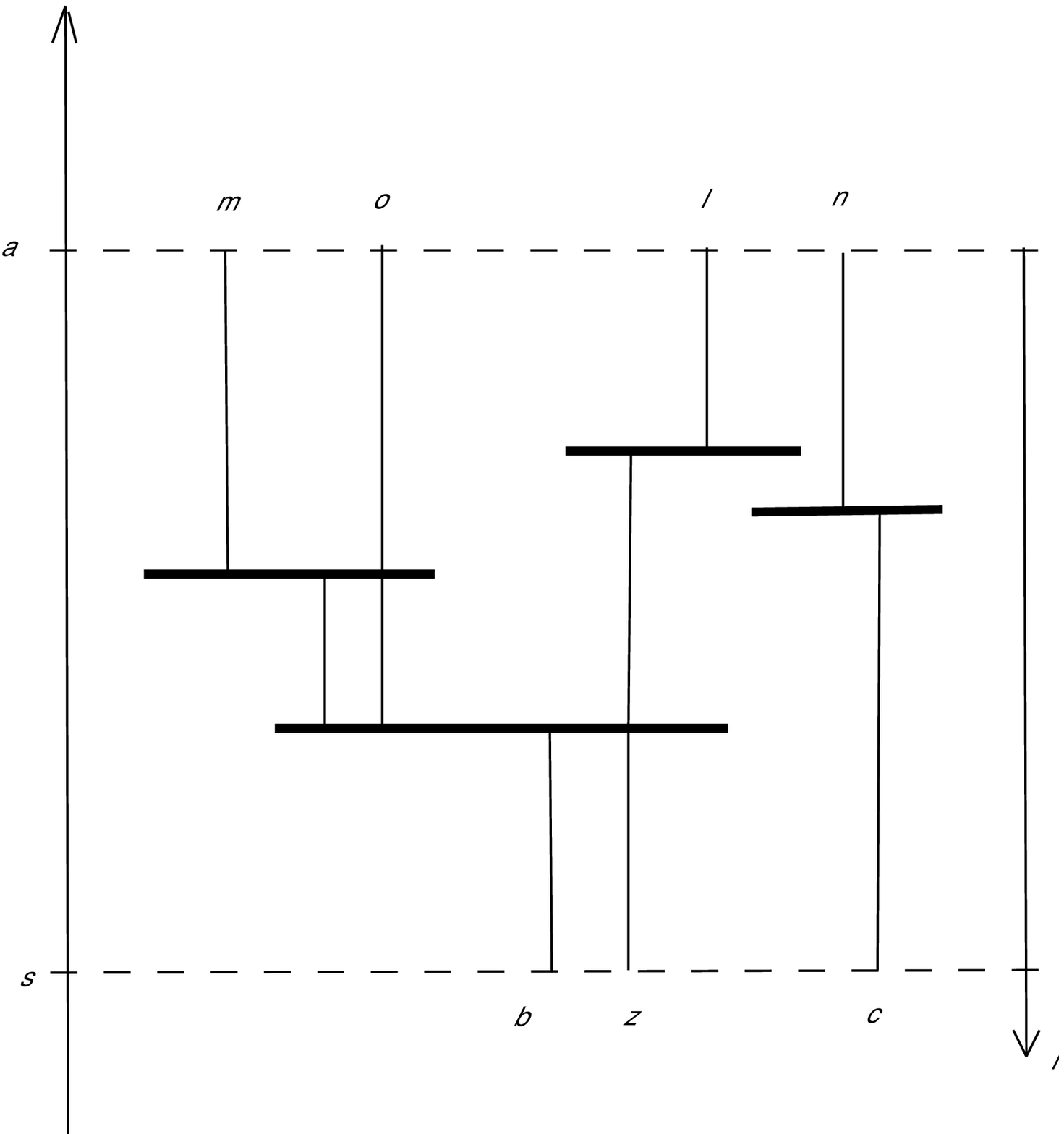}
\captionof{figure}{}
\label{ancestries}
\end{minipage}
\\\\
Condition (\ref{cond def}) guarantees that for any given $t\in \R$, for $\P$-a.e. environment $\omega$, with probability 1 no lineage in $\cA$ has an accumulation point of jumps in the time interval $[0,\infty)$. Hence, for $\mathbb P$-a.all $\omega$ we can define  $\rmPot_{\wp_k({\bf x})}$ (for all $k \in \mathbb N$ and $\lambda^{\otimes k}$-a.all $\mathbf x \in \mathbb T^k$) on the space $\mathcal{D}$ of coalescing c\`adl\`ag paths with values in $\T$.
We shall write $\rmP_{\wp_k({\bf x})}^t$ for the joint distribution on $\Omega\times \mathcal{D}$ defined by
\begin{equation}\label{joint law}
\rmP_{\wp_k({\bf x})}^t(d\omega,da):= \P(d\omega)\, \mathrm{P}_{\wp_k({\bf x})}^{\omega,t}(da).
\end{equation}
For $a \in \mathcal D$ and $h > 0$, we write $a_h$ for the restriction of $a$ to the time interval $[0,h]$ and note that $a_h$ describes a forest of coalescing paths. Indeed, if $a$ by time $h$ has coalesced to $n_h$ paths, then $a_h$ consists of $n_h$ trees $a_h^\iota$, $\iota=1,\ldots, n_h$. (Again, recall Figure \ref{ancestries}.)
\begin{rema}\label{rem: E}
Since the law $\P$ is invariant by translation in the time coordinate, for any bounded measurable function $F: \mathcal{D}\to \mathbb R$,  the quantity
$$
\int_{\Omega\times \mathcal{D}} F(a) \,\rmP_{\wp_k({\bf x})}^t(d\omega,da)
$$
is independent of $t$ and will hence be written as
$$\rmE_{\wp_k({\bf x})}[F(\mathcal A)] =
\int_{\Omega\times \mathcal{D}} F(a)\, \rmP_{\wp_k({\bf x})}(d\omega,da).
$$
(See, e.g., Corollary~\ref{co: existence}.)
\end{rema}

Another ingredient we shall need is a mutation mechanism. Let $(\cK_t)_{t\geq 0}$ be a Feller process with values in $\mathbb{K}$, defined on some probability space $(\tilde{\mathcal{D}},\tilde{\mathcal{F}},\Q)$. For every $\kappa\in \mathbb{K}$ and every genealogical tree $a$ (rooted in a single individual and) having $n$ \emph{leaves} at some time $h>0$, let us write
\begin{equation}\label{mutation}
\Q^a_\kappa\bigg[\prod_{j=1}^n g_j(\cK_h^j)\bigg],\qquad \qquad g_1,\ldots,g_n \in C(\mathbb K)
\end{equation}
to characterize the distribution of the types at the leaves when the root has type $\kappa$ and types evolve along the branches of $a$ according to the mutation process $\cK$ (we assume that this evolution occurs independently along distinct subtrees emanating from the same vertex).

The $\M_{\lambda}$-valued SLFV with mutation can now be defined (in its \emph{quenched} version) as follows.
\begin{thm1}\label{th: existence}
For $\mathbb P$-almost all $\omega$ there exists a unique $\M_{\lambda}$-valued time-inhomogeneous Hunt process $(M_t)_{t\in \R}$ whose two-parameter semigroup is characterized as follows:
For every $s \le t = s+h \in \mathbb R$, $m\in \M_{\lambda}$, $k\in \N$, $F\in C_c(\T^k)$ and $g_1,\ldots,g_k\in C(\mathbb K)$,
\begin{align}
&\rmE_{s,m}^\omega\Big[\big\la  M_t^{\otimes k}, F\otimes G_{{\bf g}}\big\ra\Big] \nonumber \\
& = \int_{\T^k} F(x_1,\ldots, x_k)\,\rmE_{\wp_k({\bf x})}^{\omega, t}\Bigg[\int_{\mathbb K^{N_h}} \prod_{\iota=1}^{N_h}\mathbb{Q}^{A_h^\iota}_{\kappa_\iota}\bigg[\prod_{j\in B^\iota_h}g_j(\cK^j_h)\bigg] \rho(\xi^1_h,d\kappa_1)\cdots \rho(\xi^{N_h}_h,d\kappa_{N_h})\Bigg] dx_1\cdots dx_k, \label{omegaduality}
\end{align}
where $\wp_k(\bf x)$ was defined in (\ref{def pn}) and $m= dx\, \rho(x,\cdot)$ as in (\ref{corr}).
\end{thm1}

Note that the right hand side of \eqref{omegaduality} is well-defined. Indeed, it follows from the above described construction of $\mathcal A$ that, conditionally on $N_h$,  the law of the ancestral  locations $\xi_h^1,\ldots,\xi_h^{N_h}$ is absolutely continuous w.r. to Lebesgue measure.

In Section~\ref{section:exist} we shall prove a slightly stronger version of Theorem 1.
Using the homogeneity of the Poisson distribution $\mathbb P$, we shall then be  able to conclude the following {\em annealed} version of  Theorem 1 (recall the notation $\rmE_{\wp_k({\bf x})}$ from Remark~\ref{rem: E}):
\begin{coro}\label{co: existence}
There exists a unique $\M_{\lambda}$-valued Hunt process $(M_t)_{t\geq 0}$ such that for every $m\in \M_{\lambda}$, $t\geq 0$, $k\in \N$, $F\in C_c(\T^k)$ and $g_1,\ldots,g_k\in C(\mathbb K)$,
\begin{align}
&\rmE_{m}\Big[\big\la  M_t^{\otimes k}, F\otimes G_{{\bf g}}\big\ra\Big] \nonumber \\
& = \int_{\T^k} F(x_1,\ldots, x_k)\,\rmE_{\wp_k({\bf x})}\Bigg[\int_{\mathbb K^{N_t}} \prod_{\iota=1}^{N_t}\mathbb{Q}^{A_t^\iota}_{\kappa_\iota}\bigg[\prod_{j\in B^\iota_t}g_j(\cK^j_t)\bigg] \rho(\xi^1_t,d\kappa_1)\cdots \rho(\xi^{N_t}_t,d\kappa_{N_t})\Bigg] dx_1\cdots dx_k. \label{duality}
\end{align}
\end{coro}
\begin{rema}
When there are no mutations (i.e., when $\mathcal{K}$ is the constant process), we recover formula (\ref{previous dual}), and hence the existence and uniqueness result from \cite{BEV2010}.
\end{rema}

Theorem~1 shows that the \emph{quenched} spatial $\Lambda$-Fleming-Viot process is a strong Markov process with c\`adl\`ag paths. In the absence of mutation, it has even stronger path properties.
\begin{lemma}\label{lem: fv}
For $\P$-a.e. environment $\omega$, the \emph{quenched} SLFV without mutation has paths of finite variation $\rmP^\omega$-a.s.
\end{lemma}
Lemma~\ref{lem: fv} is proved at the end of Section~\ref{subs: particle approx}.

\subsection{Extending the environment}
Before establishing the main results of this work, we gather here two constructions based on particular extensions of the environment. The first one will prove useful in the next sections, whereas the main interest of the second one is to relate the SLFV to Bertoin and Le Gall's flows of bridges. As discussed in the introduction, they also highlight the different layers of randomness which appear in the construction of the spatial $\Lambda$-Fleming-Viot process.

\subsubsection{Recording the locations of the parents} \label{genealogy}
Recall from the paragraph above (\ref{joint law}) that the probability measures $\rmP_{\wp_k({\bf x})}^{\omega, t}$ are defined $\P$-a.s, for each fixed $t$. As a preparation of the proof of Theorem 1, we shall specify a construction of a ``good version'' of $\rmPot_{\wp_k({\bf x})}$ which will serve simultaneously for all  $t\in \mathbb R$.

To this purpose, we mark each event $(t_i, z_i, r_i, u_i)$ with a parental location $y_i$, uniformly chosen from $\mathbb B(z_i,r_i)$, and denote the realizations of the resulting Poisson point process  by $\psi = \{(t_i, z_i, r_i, u_i, y_i): i \in I\}$. Its distribution will, by a slight abuse of notation, again be denoted by  $\mathbb P$.

For a given $\psi$, we first construct a random graph $\mathcal G^\psi$ that codes the ancestral relationships between the parental individuals living at times $t_i$ at locations $y_i$. To this purpose, let $\mathcal H := (H_{ii'})$ be a family of independent, uniformly on $[0,1]$ distributed random variables, indexed by the pairs $i\neq i'$.  For $t_{i'} < t_i$, let us define
$$G_{ii'}:= \ind_{\{H_{ii'} \ge 1-u_{i'}\}} \ind_{\{y_i\in \bB(z_{i'},r_{i'})\}}.$$
Thanks to the integrability condition \eqref{cond def},  for $\mathbb P$-almost all $\psi$ and all $i$, the random  configuration of time points $\{t_{i'}: t_{i'} < t_i, G_{ii'} =1\}$ is locally finite a.s. Thus, for all such $\psi$, there is an $\mathcal H$-measurable event of full probability on which for all $i \in I$ there is a uniquely defined index $\pi(i) \in I$ (depending on $\psi$ and $\mathcal H$) such that
$$t_{\pi(i)} = \sup\{t_{i'}:  t_{i'} < t_i, G_{ii'} =1\}.$$
\begin{minipage}{9cm}
We now decree that the most recent event which affected a parental individual with index $i$ was the event with index $\pi(i)$. In other words, the location of the ancestral lineage of the individual that is parental in the $i$-th event, when viewed backwards in time,  sits at location $y_i$ for the time span $t_i -t_{\pi(i)}$ and then jumps to $y_{\pi(i)}$. From there, it evolves as the ancestral process of the parental individual with index $\pi(i)$ (with which it has merged at time $t_i - t_{\pi(i)}$). Doing this for all $i\in I$ gives us the genealogical tree of all the parents chosen during an event, which we call $\mathcal{G}^\psi$.
\end{minipage}
\qquad
\begin{minipage}{7cm}
\psfrag{r}{$h$}
\psfrag{t}{$t$}
\psfrag{u}{$t_i$}
\psfrag{y}{$y_i$}
\psfrag{z}{$y_{\pi(i)}$}
\psfrag{s}{$t_{\pi(i)}$}
\includegraphics[width=7cm]{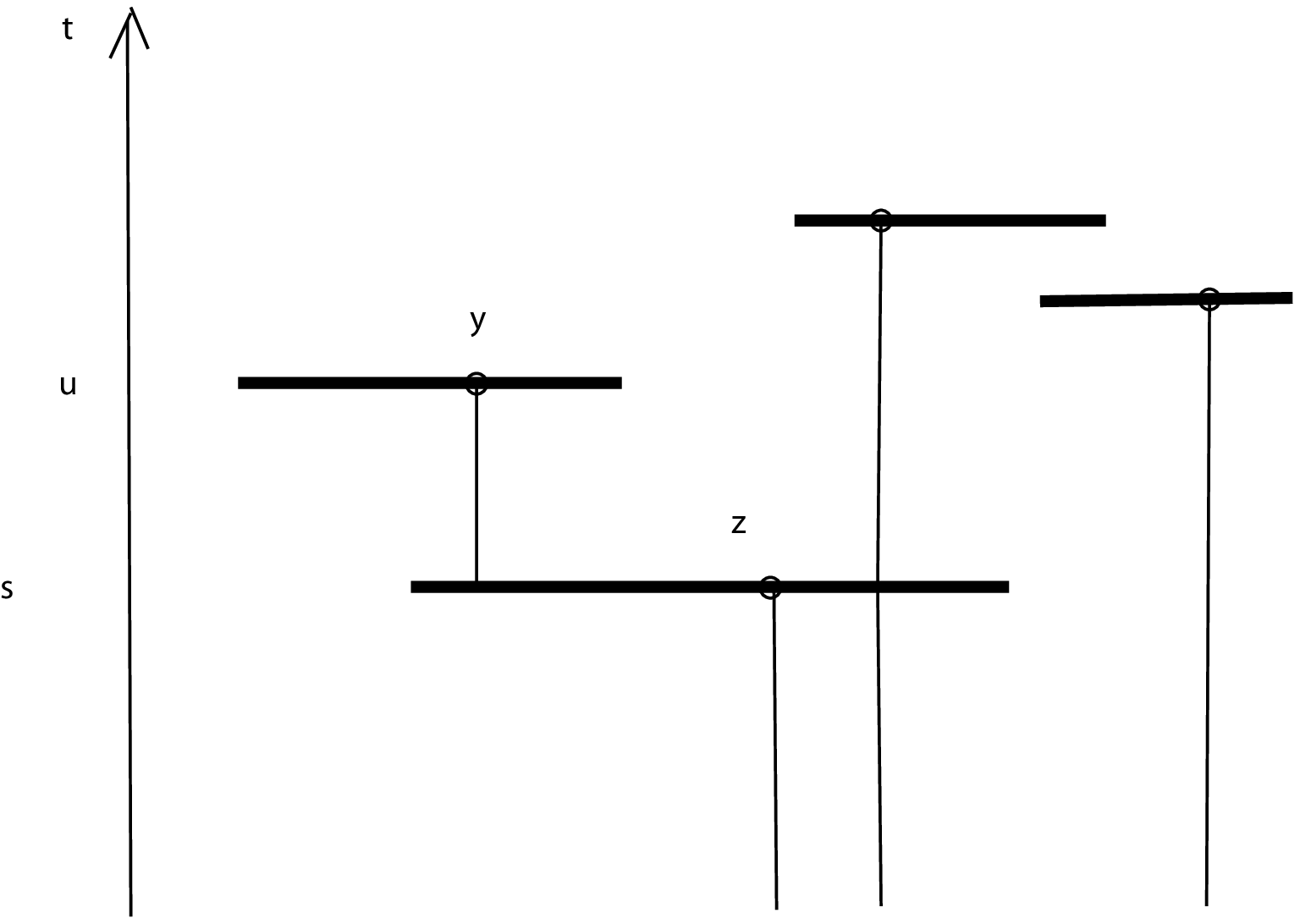}
\captionof{figure}{}
\label{parental}
\end{minipage}
\\\\
Next, let us tie in the (locations of the) ancestral lineage of an individual sampled at time $t$ at location $x$ by deciding when it joins the parental graph $\mathcal{G}^\psi$. To this purpose, we use, for given $\psi$ and $x$, a family $\mathcal H^x$ of independent, uniformly on $[0,1]$ distributed random variables $(H^x_i)_{i \in I}$ that are independent of the collection $\mathcal H$, and set
$$G_i^x:= \ind_{\{H^x_i \ge 1-u_i\}} \ind_{\{x\in \bB(z_i,r_i)\}}.$$
Again thanks to  the integrability condition \eqref{cond def}, for all $x$ and for $\mathbb P$-almost all $\psi$, the random  configuration of time points $\{t_{i}: G^x_i =1\}$ is locally finite a.s. Hence, for all $x$ and all such $\psi$, we have an $\mathcal H^x$-measurable event of full probability on which for all $t \in \mathbb R$ there is a uniquely defined index $\pi(t,x) \in I$ such that
$$t_{\pi(t,x)} = \sup\{t_{i}:  t_{i} \le t, G^x_i =1\}.$$
The location process $A^{\psi,t}_x$ of the ancestral lineage of an individual sampled at time $t$ at location $x$ is then constructed as follows. Viewed backwards in time, the lineage remains located at $x$ down to time $t_{\pi(t,x)}=:t_i$, and at that time jumps to the parental location $y_i$ in the event $e_i$. From that moment on, it has merged with the ancestral line of this parent and follows its evolution in $\mathcal{G}^\psi$.

This construction  extends to finitely or countably many individuals $\mathcal I_j$, $j \in J$, where $\mathcal I_j$ is drawn from a prescribed location $x_j$ at time $t$, and $J$ is a (finite or countably infinite) index set. For this we use independent uniformly on $[0,1]$ distributed  random variables $H^{(j)}_i$ and insert $H^{(j)}_i$ in place of $H^{x_j}_i$  to define the ancestral location process of $\mathcal I_j$.
For $\mathbf x = (x_1,\ldots, x_k)$, we then obtain the desired version of $\mathrm P_{\wp_k(\bf x)}^{\omega, t}$ by averaging the distribution of $(A^{\psi,t}_{x_1},\ldots,A^{\psi,t}_{x_k})$ over $(y_i)_{i\in I}$. Let us also note that in this way we obtain, for a Poisson point process $\Phi$ on $\mathbb T$ with Lebesgue intensity measure,  an a.s. construction of $A^{\psi,t}_{\Phi}$, the locations of the ancestral lineages of individuals sampled at the points of $\Phi$ at time $t$. This will be used at the beginning of the proof of Theorem~1'.

\subsubsection{Recording the labels of the parents} \label{flow}
This subsection, though not required for the remainder of the paper, is intended to connect the SLFV setting to the representation of generalized Fleming-Viot processes in terms of bridges by Bertoin and Le Gall \cite{BLG2003}.
To this purpose we define another extension of the environment, by assigning to each event $(t_i,z_i,r_i,u_i,y_i)$ a {\em parental label} $l_i$, independently and uniformly drawn from $[0,1]$. The realizations of the resulting Poisson point process will be denoted by $\chi = \{(t_i, z_i, r_i, u_i, y_i, l_i): i \in I\}$, and its distribution will again be denoted by $\mathbb P$.

The labels $l_i$ will be used to encode differently the additional randomness obtained through the auxiliary random variables $H_{ii'}$ and $H_i^x$ in the  construction described in the previous section. More precisely, following the ``flow of bridges'' idea of Bertoin and Le Gall \cite{BLG2003}, we shall attribute a label to each individual in the sample and use the parental labels $l_i$ to trace back the desired ancestries. As in \cite{BLG2003} we define the {\em elementary bridge} associated with $(u, l)$ as
$$
b_{l, u}(w) := (1-u) w + u\, \mathbf 1_{[l, 1]}(w), \quad 0\le w \le 1.
$$
The parental label of an individual with label $v$ that is overlapped by an event with impact $u$ and parental label $l$ is defined to be $(b_{l, u})^{-1}(v)$. In other words, under the inverse mapping $(b_{l, u})^{-1}$ the interval $[(1-u)l, (1-u)l + u]$ is mapped back to $l$, whereas the intervals $[0, (1-u)l)$ and $((1-u)l + u, 1]$ are stretched  to $[0,l)$ and $(l, 1]$, respectively.

\begin{figure}[t!]
\centering
\includegraphics[width=17cm]{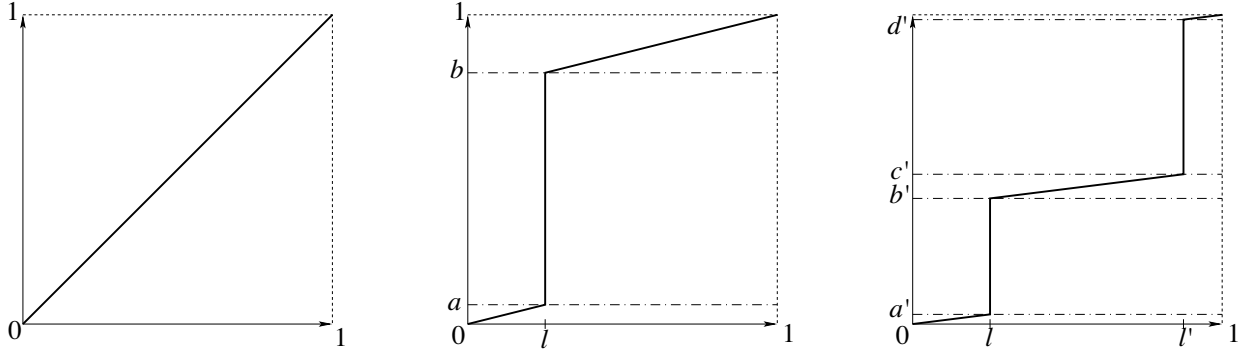}
\caption{Flow representation of the population sitting at a given site after no, 1 and 2 events. After the first event, all the individuals whose levels belong to the interval $[a,b]$ are offspring of the parent (a.s. chosen somewhere else in the area of the event) with label $l$. After the second event, this family has shrunk to the interval $[a',b']$, while the new parent with label $l'$ has given birth to all the individuals in the interval $[c',d']$. }
\label{fig: flow}
\end{figure}

For each $x\in \mathbb T$ and $s <t$, the bridge $B^x_{s,t}$ is designed below in such a way that $(B^x_{s,t})^{-1}(v)$ is the individual label of the ancestor (of the individual with label $v$ sampled at location $x$ at time $t$) that either lived at time $s$, or was parental in the most recent event before $t$ which affected the ancestral line of the individual sampled at   location $x$ at time $t$. Each event affecting $x$ between $s$ and $t$ contributes a jump to $B^x_{s,t}$ and at the same time shrinks the ``old bridge'' (see Figure~\ref{fig: flow}). Since the parental labels are a.s. distinct and since the locations of the parents are a.s. different from $x$, when looking at the population at site $x$ forwards in time, what we are trying to encode here is an immigration of jumps gradually replacing the continuous part of $B^x_{s,t}$ (which in turn represents the population at time $t$ whose ancestors at time $s$ were already living at $x$).

More concretely, the way how two elementary bridges $b_{u_1,l_1}$ and $b_{u_2,l_2}$, with $b_{u_1,l_1}$ older than $b_{u_2,l_2}$, affect $x$ is not the usual composition $b_{u_2,l_2}\circ b_{u_1,l_1}(w) = b_{u_2,l_2}(b_{u_1,l_1}(w))$, but is given by
\begin{align*}
b_{u_2,l_2}\diamond b_{u_1,l_1} (w) &:= (1-u_2)b_{u_1,l_1}(w)+u_2\mathbf 1_{[l_2, 1]}(w) \\
&= (1-u_2)(1-u_1)(w) + (1-u_2)u_1\mathbf 1_{[l_1, 1]}w + u_2\mathbf 1_{[l_2, 1]}(w), \quad 0\le w \le 1.
\end{align*}
That is, a jump of size $u_2$ is inserted at the argument $l_2$, and the old bridge is shrunk by the factor $1-u_2$ (see Figure~\ref{fig: flow}).
Doing this for all the affecting events (in the right order) gives  the bridge
\begin{align}
\label{Bst} B^x_{s,t}(w):=\left(\prod\limits_{j:s\le t_j\le t} (1-u_j)\right)w + \sum\limits_{i: s\le t_i\le t}\left(\prod\limits_{j:t_i< t_j\le t} (1-u_j)\right)u_i 1_{[l_i, 1]}(w), \quad 0\le w \le 1.
\end{align}
The quantity $p_0:= \prod\limits_{j:s\le t_j\le t} (1-u_j)$ is the fraction in the population  that is left over  at time $t$ from the population residing at location $x$ at time $s$, and the quantity $p_i := \left( \prod\limits_{j:t_i< t_j\le t} (1-u_j)\right)u_i$ is the  fraction in the population   that is left over  from the re-colonizers in the event $e_i$. In particular, in the absence of mutation (or if we only consider the family structure of the population at a given site), we see that the type (or family size) distribution conditioned on the extended environment $\chi$ is a deterministic function of its initial value and of $\chi$. Furthermore, for a random variable $V$ that is uniformly distributed on $[0,1]$, the probability that the pre-image $(B^x_{s,t})^{-1}(V)$ equals $l_i$ is  $p_i$, and the probability that it is a continuity point of $B^x_{s,t}$ is $p_0$.

Based on this last point, let us show that the ancestry of a sample of individuals can be reconstructed by a deterministic procedure involving only $\chi$ and the distribution of the level $V_j$ of the individuals that we consider. Indeed, on the event $(B^x_{s,t})^{-1}(V)= l_i$ for some $i$ indexing an event, we set $T(x,V):= t_i$, and proceed back from $(t_i, y_i)$ and from the individual level $l_i$, this time using those of the extended events in $\chi$ that affect the location $y_i$.

The integrability condition implies that for all $x$ and $\mathbb P$-a.a. $\chi$
$$\sum_{i: s\le t_i \le t,\, |z_i-x|\le r_i} u_i < \infty.$$
Therefore, for all $v\in [0,1]$, the path $h \mapsto (B^x_{t-h,t})^{-1}(v)$, $0\le h \le t$, is a.s. of finite variation (though with jumps that generically will occur densely in time). Let us also note that for each $h \in [0,t]$,  conditionally on $h> t- T(x,V)$, the random variable $(B^x_{t-h,t})^{-1}(V)$ is uniformly distributed on $[0,1]$ when randomized over $\chi$.

We can carry out the same construction for a sequence $(x_j, V_j)$ with independent $V_1, V_2, \ldots$, and the locations $x_j$ not necessarily distinct. In this way we obtain a $(\chi, V_1, V_2, \ldots)$-measurable configuration of coalescing ancestral lines starting in $(x_j, V_j)_{j \ge 1}$ at time $ t$, whose $x$-component, when randomized over $\chi$, has the distribution $\mathrm P_{\wp_k(\bf x)}^{t}$.

\section{The SLFV conditioned on the environment}\label{section:exist}

In this section, we prove a slightly stronger version of Theorem~1, which will establish the latter \emph{a fortiori}: for $\mathbb P$-almost all $\omega \in \Omega$, we characterize a two-parameter semigroup $Q_{s,t}^\omega (m, dm')$ of transition probabilities on $\mathcal M_\lambda$, where  $Q_{s,t}^\omega (m, \cdot)$ is the distribution of a  random element $\fm_{s,t}(m)$ in $\mathcal M_\lambda$, which is the result of the transformation of $m$ due to the mutations and the events in $\omega$ occurring in the time interval $[s,t]$. Notice that the semigroup $Q_{s,t}^\omega$ is not time-homogeneous, since the sequence of reproduction events encoded in the environment $\omega$ is not time-homogeneous.
\begin{thm1prime}\label{prop: flow}
For $\mathbb P$-almost every $\omega\in {\Omega}$, there exists a unique collection $Q_{s,t}^\omega (m, dm')$ of transition probabilities on $\mathcal M_\lambda$ such that for every $s \le t \in \R$,  $k\in \N$, $F\in C_c(\T^k)$, $g_1,\ldots,g_k\in C(\mathbb K)$, and $h:= t-s$,
\begin{align}
& \int Q_{s,t}^\omega (m, dm')\big\la m'^{\otimes k}, F\otimes G_{{\bf g}}\big\ra \label{omega duality} \\
& = \int_{\T^k} F(x_1,\ldots, x_k)\,\rmE^{\omega,t}_{\wp_k({\bf x})}\Bigg[\int_{\mathbb K^{N_h}} \prod_{\iota=1}^{N_h}\mathbb{Q}^{A_h^\iota}_{\kappa_\iota}\bigg[\prod_{j\in B^\iota_h}g_j(\cK^j_h)\bigg] \rho(\xi^1_h,d\kappa_1)\cdots \rho(\xi^{N_h}_h,d\kappa_{N_h})\Bigg] dx_1\cdots dx_k, \nonumber
\end{align}
where here again $m(dx,d\kappa)= dx\, \rho(x,d\kappa)$.
\smallskip
\noindent Furthermore, $\P$-a.s. we have

\noindent \item[(i)]  for all $s\leq t$, $Q_{s,t}^\omega (m, \cdot )$ is a continuous function of $m\in \M_{\lambda}$,

\noindent \item[(ii)] for all $s\in \R$, $t \to Q_{s,t}^\omega, \, t\ge s$, is a strongly continuous operator.

\noindent \item[(iii)] $(Q_{s,t}^\omega)$ is a two-parameter semigroup.
More precisely, for all $q < s<t$, and  for every $F$ and $g_1,\ldots,g_k$ as before we have
\begin{align}
\int Q_{q,s}^\omega (m, dm' )\int Q_{s,t}^\omega (m', dm'')\big\la  m''^{\otimes k}, F\otimes G_{{\bf g}}\big\ra   = \int Q_{q,t}^\omega (m, dm'' )\big\la m''^{\otimes k}, F\otimes G_{{\bf g}}\big\ra. \label{concat}
\end{align}
\end{thm1prime}

{\rema \label{dua}
Stated in words, relation \eqref{omega duality} reads as follows: Let $\mathfrak M_{t-h,t}(m)$ be a random measure with distribution $Q_{s,t}^\omega(m,\cdot)$. Then a random $k$-sample drawn at locations $x_1,\ldots, x_k$ from $\mathfrak M_{t-h,t}(m)$ arises in three steps:    (1) Take, under $\rmP^{\omega,t}$,  the random ancestry $\mathcal A$ started in $\wp_k({\bf x})$, (2) sample the types $\kappa_1, \ldots, \kappa_{N_h}$ drawn at locations $\xi_h^1, \ldots, \xi_h^{N_h}$ from $m$, and (3) run the mutation process  starting from the types $\kappa_1, \ldots, \kappa_{N_h}$ along the genealogy $\mathcal A$ forwards in time, starting at time $s=t-h$, and up to time $t$. The types sampled at locations $x_1,\ldots, x_k$ at time $t$ then result from the forest-indexed mutation process that is run for $h$ units of time.}
\\\\

Before we prove Theorem 1', let us show that Corollary~\ref{co: existence} is a straightforward consequence of it.

\medskip
\noindent\textbf{Proof of Corollary~\ref{co: existence}.} For any initial condition $m\in \M_{\lambda}$, let us set $Q_t(m,\cdot)= \int \mathbb P(d\omega) Q^\omega_{0,t}(m,\cdot)$.
Using (\ref{omega duality}), \eqref{joint law} and Fubini's theorem,
 we obtain that for every $k\in \N$, $F\in C_c(\T^k)$ and $g_1,\ldots,g_k\in C(\mathbb K)$,
\begin{align}
&\int Q_{t} (m, dm')\big\la m'^{\otimes k}, F\otimes G_{{\bf g}}\big\ra \nonumber \\
& = \E\Bigg[\int_{\T^k} F(x_1,\ldots, x_k)\,\rmE^{\omega,t}_{\wp_k({\bf x})}\Bigg[\int_{\mathbb K^{N_t}} \prod_{\iota=1}^{N_t}\mathbb{Q}^{\cA_t^\iota}_{\kappa_\iota}\bigg[\prod_{j\in B^\iota_t}g_j(\cK^j_t)\bigg] \rho(\xi^1_t,d\kappa_1)\cdots \rho(\xi^{N_t}_t,d\kappa_{N_t})\Bigg] dx_1\cdots dx_k\Bigg], \nonumber \\
& = \int_{\T^k} F(x_1,\ldots, x_k)\,\rmE_{\wp_k({\bf x})}\Bigg[\int_{\mathbb K^{N_t}} \prod_{\iota=1}^{N_t}\mathbb{Q}^{\cA_t^\iota}_{\kappa_\iota}\bigg[\prod_{j\in B^\iota_t}g_j(\cK^j_t)\bigg] \rho(\xi^1_t,d\kappa_1)\cdots \rho(\xi^{N_t}_t,d\kappa_{N_t})\Bigg] dx_1\cdots dx_k\Bigg]. \label{randomized}
\end{align}
Hence (\ref{duality}) holds and the existence of a suitable collection $\{M_t,\, t\geq 0\}$ of random variables is proven. Next, Lemma~\ref{lem:topo}$(c)$ shows that the duality relations (\ref{duality}) are sufficient to guarantee the uniqueness in law of each $M_t$.

There remains only to show that $\{Q_t,\, t\geq 0\}$ is a semigroup. Indeed, once this has be established,  Theorem 1' $(i)$ and $(ii)$ show that $\{Q_t,\, t\geq 0\}$ has the Feller property, and hence defines a Hunt process. Now for $0\le s\le t$ let us observe that,  because of \eqref{concat}
$$Q_{t}(m,\cdot) = \int \mathbb P(d\omega) Q^\omega_{0,t}(m,\cdot) =  \int \mathbb P(d\omega) \int_{\mathcal M_\lambda}Q^\omega_{0,s}(m,dm')Q^\omega_{s,t}(m',\cdot) .$$
Since $\omega \mapsto Q^\omega_{0,s}$ is measurable with respect to the events between times $0$ and $s$, $\omega \mapsto Q^\omega_{s,t}$ is   measurable with respect to the events between times $s$ and $t$, and $\mathbb P$-a.s. no event happens precisely at time $s$, then due to the independence property of the Poisson distribution $\mathbb P$ and by Fubini's theorem we can rewrite the r.h.s. as
\begin{align}
\label{omegawise}
\int_\Omega \mathbb P(d\omega)\int_{\mathcal M_\lambda}Q^\omega_{0,s}(m,dm') \int_\Omega \mathbb P(d\omega')Q^{\omega'}_{s,t}(m',\cdot).
\end{align}
Since the environment distribution $\mathbb P$ (and together with it also the distribution of the genealogies randomized over the environment)  is invariant under time-shift, we can write that $ \int_\Omega \mathbb P(d\omega')Q^{\omega'}_{s,t}(m',\cdot) =  \int_\Omega \mathbb P(d\omega')Q^{\omega'}_{0,t-s}(m',\cdot).$ Hence \eqref{omegawise} equals $\int_{\mathcal M_\lambda}Q_s(m,dm')Q_{t-s}(m',\cdot)$, which completes the proof of the semigroup property. \hfill $\Box$

\medskip
\noindent\textbf{Proof of Theorem 1'.}
\emph{Existence of $Q_{s,t}^\omega(m)$.} In a first step,  we consider a Poisson point measure $\Phi$ on $\mathbb T$ with intensity measure $\lambda$. We write $\Phi = \sum_{j\in J} \delta_{x_{j}}$ and construct, for a given environment $\psi = (t_i,z_i,r_i,u_i,y_i)_{i\in I}$, the coalescing ancestral location processes $(A^{\psi,t}_{x_{j}})_{j\in J}$ as described in Section~\ref{genealogy}. We run this processes for $h=t-s$ units of time into the past and sample a type at each location $\xi_h^\iota$ of an ancestor according to the type distribution $\rho(\xi_h^\iota, \cdot)$. Finally, we run the mutation process described in the paragraph around (\ref{mutation}) along each of the ancestral trees we obtained, independently of each other. In other words, assuming the type $\kappa_\iota$ was attributed to the $\iota$-th ancestor, we use that allele as a starting point of the evolution along the genealogical tree rooted at this ancestor, and the result of the mutation process at time $h$ gives us the type of each of the points $x_j\in\Phi$ such that $j$ belongs to the block $B^\iota_h$. In this way, we obtain a point measure $\cN$ on $\T\times \mathbb K$, whose projection onto the geographical space forms a Poisson point process whilst the types assigned to the different points are correlated through the genealogy (and over the geography).

In the second step, we perform the just described procedure for countably many i.i.d. copies $\Phi_1, \Phi_2,\ldots$ of $\Phi$, resulting in a sequence $\cN_1,\cN_2, \ldots$ of random point measures on $\T\times \mathbb K$. A simple but crucial observation is that the construction of the point measures $\cN_n$ is equivalent in law to the following construction:
\begin{itemize}
\item[$(i)$] Enrich the environment with the family $\mathcal{H}$ and construct the parental skeleton $\mathcal{G}^\psi$ as in Section~\ref{genealogy}.
\item[$(i)$] Allocate a type to each of the parents by running the mutation process along $\mathcal{G}^\psi$, from time $s$ on. This gives us a labelled skeleton $\mathcal{S}$.
\item[$(ii)$] For every $x\in \sum_n \Phi_n$, if the ancestral location process back from $x$ coalesces with $\mathcal{G}^\psi$ (say, at time $\tau_x$ back into the past), then start an independent mutation process from the type at this point of the skeleton and run it for the time span $\tau_x$. If the ancestral process does not coalesce with $\mathcal{G}^\psi$ (and thus remains at $x$), then sample a type from the measure $\rho(x,\cdot)$. In both cases, allocate the resulting type to the `individual' located at $x$ in the point process.
\end{itemize}
This description shows that, conditionally on $\mathcal{S}$, the point measures $\mathcal{N}_1,\mathcal{N}_2,\ldots$ are independent and identically distributed. From this it is not difficult to conclude that $\frac 1n\cN^n := \frac 1n(\cN_1+\cdots +\cN_n)$  converges a.s. in the vague topology towards a (random) limit $\fm_{s,t}$ as $n\to \infty$. Indeed, let us write $\cB_c(\T\times \mathbb K)$ for the set of all bounded functions with compact support on $\T\times \mathbb K$. By the law of large numbers (applied conditionally under $\psi$ and $\mathcal S$),
for every $f\in \cB_c(\T\times \mathbb K)$, the sequence $(\la \frac{1}{n}\cN^n,f\ra)_{n\in \N}$ converges $\rmP^{\omega,t}$-a.s. and in $L^1$ to a random variable $L(f)$. We thus put,  for every $f\in \cB_c(\T\times \mathbb K)$,
\begin{equation}\label{def fm}
\la \fm_{s,t}, f\ra := L(f).
\end{equation}
Up to restricting this definition to a countable set of functions $f$ of the form $\ind_{A_j}$, where the compact sets $\{A_j,\, j\in \N\}$ form a basis of $\T\times \mathbb K$, we can conclude that $\rmP^{\omega,t}$-a.s., $\fm_{s,t}$ is a (random) nonnegative Radon measure and $\frac{1}{n}\cN^n$ converges a.s. in the vague topology towards $\fm_{s,t}$ as $n\to \infty$. On the other hand, again by the law of large numbers, for each $F \in C_c(\mathbb T)$ and $f:= F \otimes \ind_\mathbb K$,  we have that $\la \frac{1}{n}\cN^n, f\ra \to \la \lambda, F\ra$ a.s. as $n\to \infty$. Hence, with probability one $\fm_{s,t}$ belongs  to $\mathcal M_\lambda$.

Defining $Q_{s,t}^\omega (m, .)$ as the law of $\fm_{s,t}$, it remains to prove (\ref{omega duality}). For this we proceed by proving the following claim:

\noindent For every $k\in \N$, $F\in C_c(\T^k)$, $g_1,\ldots,g_k\in C(\mathbb K)$, we have
\begin{align}
\lim_{n\rightarrow \infty} &\mathrm{E}^{\omega,t}\Big[\big\la \big(n^{-1}\cN^n\big)^{\otimes k}, F\otimes G_{{\bf g}}\big\ra\Big]  \label{asdual}\\
& = \int_{\T^k} F(x_1,\ldots, x_k)\,\rmE^{\omega,t}_{\wp_k({\bf x})}\Bigg[\int_{\mathbb K^{N_h}} \prod_{\iota=1}^{N_h}\mathbb{Q}^{\cA_h^\iota}_{\kappa_\iota}\bigg[\prod_{j\in B^\iota_h}g_j(\cK^j_h)\bigg] \rho(\xi^1_h,d\kappa_1)\cdots \rho(\xi^{N_h}_h,d\kappa_{N_h})\Bigg] dx_1\cdots dx_k. \nonumber
\end{align}
Indeed, writing $\cN^n = \sum_\gamma \delta_{(X_\gamma,K_\gamma)}$ we can compute
\begin{align*}
\rmE^{\omega,t}\bigg[\bigg\la \frac{1}{n^k}\big(\cN^n\big)^{\otimes k},\, F\otimes G_{\bf g}\bigg\ra\bigg] & = \frac{1}{n^k}\,\rmE^{\omega,t}\bigg[\sum_{\gamma_1,\ldots,\gamma_k}F\big(X_{\gamma_1},\ldots,X_{\gamma_k}\big)\prod_{\jmath=1}^k g_\jmath\big(K_{\gamma_\jmath}\big)\bigg]\\
&=\frac{1}{n^k}\, \rmE^{\omega,t}\bigg[\sum_{\gamma_1\neq \ldots\neq \gamma_k}F\big(X_{\gamma_1},\ldots,X_{\gamma_k}\big)\prod_{\jmath=1}^k g_\jmath\big(K_{\gamma_\jmath}\big)\bigg] + \mathcal{O}\Big(\frac{1}{n}\Big)
\end{align*}
as $n\rightarrow \infty$ (the sum is set to be zero if $\cN^n$ has less than $k$ points, which happens only if $\mathrm{Vol}(\T)<\infty$ and has a probability tending to $0$ as $n\rightarrow \infty$ in this case). Writing $\Phi^n:= \Phi_1+\cdots + \Phi_n$, we have
\begin{align*}
\frac{1}{n^k}\, \rmE^{\omega,t}\bigg[&\sum_{\gamma_1\neq \ldots \neq \gamma_k}F\big(X_{\gamma_1},\ldots,X_{\gamma_k}\big)\prod_{\jmath=1}^k g_\jmath\big(K_{\gamma_\jmath}\big)\bigg] \\
& = \frac{1}{n^k}\, \rmE^{\omega,t}\bigg[\rmE^{\omega,t}\bigg[\sum_{\gamma_1\neq \ldots\neq \gamma_k}F\big(X_{\gamma_1},\ldots,X_{\gamma_k}\big)\prod_{\jmath=1}^k g_\jmath\big(K_{\gamma_\jmath}\big)\, \Big|\, \Phi^n\bigg]\bigg]\\
& = \frac{1}{n^k}\, \rmE^{\omega,t}\bigg[\sum_{\gamma_1\neq \ldots\neq \gamma_k}F\big(X_{\gamma_1},\ldots,X_{\gamma_k}\big)\rmE^{\omega,t}\bigg[\prod_{\jmath=1}^k g_\jmath\big(K_{\gamma_\jmath}\big)\, \Big|\, \Phi^n\bigg]\bigg].
\end{align*}
But by construction, conditionally on the distinct values of $X_{\gamma_1},\ldots,X_{\gamma_k}$ we have, again with $h:= t-s$,
\begin{align*}
\rmE^{\omega,t} \, \bigg[\prod_{\jmath=1}^k g_\jmath\big(K_{\gamma_\jmath}\big)\, \Big|\,\Phi^n\bigg]& =  \rmE^{\omega,t}_{\wp_k({\bf X})}\Bigg[\int_{\mathbb K^{N_h}} \prod_{\iota=1}^{N_h}\mathbb{Q}^{\cA_h^\iota}_{\kappa_\iota}\bigg[\prod_{\jmath\in B^\iota_h}g_\jmath(\cK^j_h)\bigg] \rho(\xi^1_h,d\kappa_1)\cdots \rho(\xi^{N_h}_h,d\kappa_{N_h})\Bigg]  \\& =: \Psi_g({\bf X}).
\end{align*}
Thus, using the $k$-th moment formula for the Poisson point measure $\Phi^n= \sum_\gamma \delta_{X_\gamma}$ we arrive at
\begin{align*}
\rmE^{\omega,t}\, \bigg[\bigg\la \frac{1}{n^k}\big(\cN^n\big)^{\otimes k},\, F\otimes G_{\bf g}\bigg\ra\bigg]& = \frac{1}{n^k}\, \rmE^{\omega,t}\bigg[\sum_{\gamma_1\neq \ldots \neq \gamma_k}F\big(X_{\gamma_1},\ldots,X_{\gamma_k}\big)\Psi_g\big(X_{\gamma_1},\ldots,X_{\gamma_k}\big)\bigg] +\mathcal{O}\Big(\frac{1}{n}\Big) \\
& = \int_{\T^k}F(x_1,\ldots,x_k)\Psi_g(x_1,\ldots,x_k)\, dx_1\cdots dx_k + \mathcal{O}\Big(\frac{1}{n}\Big),
\end{align*}
which proves the claimed equality \eqref{asdual}.

Denoting the (compact) support of the function $F$ by $\mathrm{Supp}(F)$ and using the fact that the number of points of $\Phi^n$ within the support of $F$ is a Poisson random variable with parameter $n\mathrm{Vol}(\mathrm{Supp}(F))$, we can write
\begin{align*}
\mathrm{E}^{\omega,t}\Big[\big|\big\la\big(n^{-1}\cN^n\big)^{\otimes k},F\otimes G_{{\bf g}} \big\ra\big|\Big] & \leq \bigg\{\|F\|_{\infty}\prod_{\jmath=1}^k\|g_\jmath\|_{\infty}\bigg\}\, \mathrm{E}^{\omega,t}\Big[\big\la\big(n^{-1}\cN^n\big)^{\otimes k},\ind_{\mathrm{Supp}(F)} \big\ra\Big]\nonumber \\
& \leq \bigg\{\|F\|_{\infty}\prod_{\jmath=1}^k\|g_\jmath\|_{\infty}\bigg\}\, \frac{C n^k\mathrm{Vol}(\mathrm{Supp}(F))^k}{n^k} \\
& = C \bigg\{\|F\|_{\infty}\prod_{\jmath=1}^k\|g_\jmath\|_{\infty}\bigg\}\, \mathrm{Vol}(\mathrm{Supp}(F))^k
\end{align*}
for a constant $C$ independent of all other parameters. We can therefore use dominated convergence, together with the fact that $\frac{1}{n}\cN^n$ converges vaguely to $\fm_{s,t}$ with $\rmP^{\omega,t}$-probability $1$, to conclude that
$$
\lim_{n\rightarrow \infty} \mathrm{E}^{\omega,t}\Big[\big\la \big(n^{-1}\cN^n\big)^{\otimes k}, F\otimes G_{{\bf g}}\big\ra\Big]  = \rmE^{\omega,t}\Big[\big\la \fm_{s,t}^{\otimes k},F\otimes G_{{\bf g}} \big\ra\Big].
$$
Combining the above with \eqref{asdual} yields (\ref{omega duality}), and the proof of existence is complete.

\medskip
\emph{Uniqueness.} By Lemma~\ref{lem:topo}$(c)$, the equalities (\ref{omega duality}) for all functions of the form $I_k$ are sufficient to ensure that there is at most one distribution on $\M_{\lambda}$ which satisfies them for any fixed $\omega$, $m$, $s$, $t$. Hence, uniqueness holds.

\medskip
\emph{Continuity with respect to $m$.} This is a direct consequence of (\ref{omega duality}). Indeed, knowing that $k$ lineages are sampled from the compact support of $F$, one can truncate the distribution of $(\xi^1_t,\ldots,\xi^{N_t}_t)$ uniformly in $(x_1,\ldots,x_k)\in \mathrm{Supp}(F)$ and turn the truncated density into a continuous function with compact support so that, up to an arbitrarily small error term, the expectation in the right-hand side of (\ref{omega duality}) is an integral w.r. to $m^{\otimes N_h}$ of some continuous and compactly supported function. Assume now that a sequence $(m_n)_{n\geq 1}$ of elements of $\M_{\lambda}$ converges vaguely towards $m$ (and so $m_n^{\otimes k}$ converges vaguely towards $m^{\otimes k}$ for every $k$). By the definition of vague convergence, for every $(x_1,\ldots,x_k)\in \mathrm{Supp}(F)$ the expectation in the right-hand side of (\ref{omega duality}) written with $m_n$ converges to the same expectation with $m$ as the initial type distribution. Using dominated convergence (and working `up to an arbitrarily small correction term')  we can conclude that
$$
\lim_{n\to \infty} \int Q_{s,t}^\omega (m_n, dm' )\big\la  m'^{\otimes k}, F\otimes G_{{\bf g}}\big\ra =\int Q_{s,t}^\omega (m, dm')\big\la m'^{\otimes k}, F\otimes G_{{\bf g}}\big\ra.
$$
By Lemma~\ref{lem:topo}$(c)$, this guarantees that $Q_{s,t}^\omega (m_n,\cdot)$ converges weakly to $Q_{s,t}^\omega (m, \cdot)$ as $n\rightarrow \infty$.

\medskip
\emph{Strong continuity at $s$.} We follow the construction in Subsection \ref{genealogy}. Here the good $\omega$'s are the ones for which the parental skeleton exists (with a locally finite jump intensity along all its lineages) {\em and} for which for $\lambda$-almost all $x$ the set of jump points $\{t_i: G_i^x = 1\}$ is locally finite on the time axis.

By Lemma~\ref{lem:topo}$(c)$, all we need to show is that for all $\omega$ with the just described property we have
\begin{equation}\label{aim}
\lim_{\e\rightarrow 0+}  \int Q_{s,s+\varepsilon}^\omega (m, dm' )\big\la  m'^{\otimes k}, F\otimes G_{{\bf g}}\big\ra=  \big\la  m^{\otimes k}, F\otimes G_{{\bf g}}\big\ra
\end{equation}
for every $m\in \M_\lambda$, $k\in \N$, $F\in C_c(\T^k)$ and $g_1,\ldots,g_k\in C(\mathbb K)$.

Here again, the key tool is the set of equations (\ref{omega duality}).
We can couple the ${\rm P}^{\omega,s+\varepsilon}_{\wp_k(\mathbf x)}$ just as we did in Subsection \ref{genealogy} to see that for Lebesgue-almost all $\mathbf x = (x_1,\ldots, x_k)\in \mathrm{Supp}(F)$,
$$
\lim_{\e\rightarrow 0} \rmP^{\omega,s+\e}_{\wp_k(\mathbf{x})}\big[\xi \hbox{ jumps during }(s,s+\e]\big]
= 0.
$$
Together with the Feller property of the mutation process, we can thus conclude that  for Lebesgue-a.e. $\mathbf{x}\in \mathrm{Supp}(F)$
$$
\lim_{\e\rightarrow 0} \rmE^{\omega,s+\e}_{\wp_k({\bf x})}\Bigg[\prod_{\iota=1}^{N_\e}\mathbb{Q}^{A_\e^\iota}_{\kappa_\iota}\bigg[\prod_{j\in B^\iota_\e}g_j(\cK^j_\e)\bigg]\bigg] = \prod_{j=1}^k g_j(\kappa_j).
$$
This result being independent of the measure $m$ taken to attribute the types $\kappa_1,\ldots,\kappa_k$, (\ref{omega duality}) and dominated convergence give us that (\ref{aim}) is satisfied for every $m\in \M_\lambda$, $\P$-a.s.

\medskip
\emph{Flow property.} Let us give a conceptual proof, based on Remark \ref{dua} (which is stated below Theorem 1'). For brevity, we write $\mathfrak M^{\omega}_{q,s}(m)$ for a random element in $\M_\lambda$ that has distribution $Q^\omega_{q,s}(m,\cdot)$. With that notation, we have to show that $\mathfrak M^{\omega}_{s,t}(\mathfrak M^{\omega}_{q,s}(m)) \stackrel d=\mathfrak M^{\omega}_{q,t}(m)$ for $\mathbb P$-almost all $\omega$. By Lemma~\ref{lem:topo}, for this it is enough to show that  the type distribution of a random $k$-sample drawn at locations $x_1,\ldots, x_k$ coincide for both random measures. We start with analysing these for $\mathfrak M^{\omega}_{s,t}(\mathfrak M^{\omega}_{q,s}(m))$. In order to obtain the random types at the locations $\xi_{t-s}^1, \ldots, \xi_{t-s}^{N_{t-s}}$  as required in step (2) of Remark \ref{dua}, we have to go further down in the ancestry, now starting from the partition $\{(\{1\},\xi_{t-s}^1), \ldots, (\{N_{t-s}\},\xi_{t-s}^{N_{t-s}})\}$. By the Markov property of $\mathcal A$ under $\rmPot_{\wp_k({\bf x})}$ (see Section~\ref{genealogy}) , this amounts to running $\mathcal A$ for a total time $t-q$. Consequently, we have to sample at the locations $\xi_{t-q}^1, \ldots, \xi_{t-q}^{N_{t-q}}$ from the measure $m$, and from the resulting types run the mutation process forwards in time between times $q$ and $s$ in order to obtain the types at time $s$. These are used as the input in step (3) of Remark \ref{dua}, that is, as the initial conditions for another go of the mutation process between times $s$ and $t$. By the Markov property of the mutation process $\mathcal K$, this amounts to running $\mathcal K$ along the genealogy $\mathcal A$ between times $q$ and $t$. In total, we have thus arrived at the types of a random $k$-sample drawn at locations $x_1,\ldots, x_k$ from $\mathfrak M^{\omega}_{q,t}(m)$. Notice that this scheme being independent of the measure $m$, it works simultaneously for all $m$'s.

The proof of Theorem 1' is now complete. \hfill $\Box$

\section{A look-down representation of the SLFV-process}\label{section:look-down}
In the previous section we constructed the two-parameter semigroup $(Q_{s,t}^\omega)$ given the configuration $\omega$ of events. This was done, first for fixed times $s < t$, on top of the random genealogy of a sample whose locations had a Poisson$(c\lambda)$-distribution (with $c \to \infty$). One ingredient for constructing this random genealogy was the process of parental locations in the events, which was independent of the sampling locations. In this section, we shall again work with a Poisson system of sampling locations ``with infinite spatial density'' but now this system will evolve in time, with the locations of the sampling and the parental locations in the events being coupled. That is, we shall always choose the parent among the individuals of our Poisson system.

Let us start with an informal description.
At time $t=0$ we start with a Poisson configuration of particles on $\mathbb T \times [0,\infty)$ with intensity measure $\lambda(dx) \otimes d\ell$. The first component is the particle's {\em location}, the second will be called the particle's {\em level}. While the levels stay fixed in time, the locations $\zeta_j(t)$, given $\omega =(t_i,z_i,r_i,u_i)_{i\in I} \in \Omega$, perform independent jump processes: at each time $t_i$ such that $||\zeta_j(t_i-) - z_i|| \le r_i$, the particle at level $\ell_j$ tosses a coin with success probability $u_i$, independently of everything else. If this coin comes up with ``success'', the particle jumps to a location $\zeta_j(t_i)$ that is chosen uniformly from $\mathbb B(z_i,r_i)$, again independently of everything else.

Because of the integrability condition \eqref{cond def}, the jump times of every single particle are a.s. locally finite, and because of the Poisson colouring theorem, the process $(\zeta_j(t), \ell_j)_j$ remains Poisson,   with invariant intensity measure $\lambda(dx) \otimes d\ell$ (as we shall see below). Among all the particles that are affected by the event at time $t_i$ there will therefore be a.s. one whose level is smallest, say $\ell_{j(i)}$. The key idea is now to declare the   location $\zeta_{j(i)}(t_i-)$ as the parental location in the event at time $t_i$, and to decree that all those $(t_i, \zeta_j(t_i))$ for which $\zeta_j$ jumped at time $t_i$ are   {\em children} of $(t_i-, \zeta_{j(i)}(t_i-))$. In this way, a genealogy is filled into  the space-time point configuration $(t_i, \zeta_j(t_i))$ in a {\em look-down} manner: certain ones of the particles at higher levels {\em look down} at a particle  at lower level  and copy its type.

Let us now proceed with a more formal definition of the particle system in order to prove a look-down representation of the SLFV as stated in Theorem 2. We shall discuss possible generalizations of this construction in Section~\ref{subs: generalizations}.

\subsection{A particle system and its look-down genealogy}\label{subs: ld}

Let us fix again an environment $\omega\in \Omega$. Let us define the \emph{forwards-in-time} motion $(\zeta_t)_{t\geq 0}$ of a single individual in the environment $\omega$ as follows: whenever $\zeta$ lies within the range of an event $(t_i,z_i,r_i,u_i)\in \omega$, the process at time $t_i$ does nothing with probability $1-u_i$, or jumps to a new location uniformly distributed over $\mathbb B(z_i,r_i)$ with probability $u_i$. We write $\rmP^{\omega}_x$ for the probability measure on $\mathcal{D}$ under which $\zeta$ starts at $x\in \T$ (recall that $\mathcal{D}$ is the space of coalescing c\`adl\`ag paths with values in $\T$, defined just above (\ref{joint law})).

Let us observe that Lebesgue measure is reversible for the evolution of $\zeta$. Indeed, writing $K_i(x,dy)$ for the transition kernel of $\zeta$ during the event of $\omega$ labelled by $i$, we have
\begin{equation}\label{kernel}
K_i(x,dy) = \ind_{\{x\notin \mathbb B(z_i,r_i)\}}\, \delta_x + \ind_{\{x\in \mathbb B(z_i,r_i)\}}\, \Big\{(1-u_i)\delta_x + \frac{u_i}{\mathrm{Vol}(\mathbb B(z_i,r_i))}\, dy\Big|_{\mathbb  B(z_i,r_i)}\Big\},
\end{equation}
so that it is easy to check that
$$
dx \, K_i(x,dy) = dy\, K_i(y,dx), \quad x,y\in \T.
$$
In particular, Lebesgue measure is conserved by the flow through the countably many events of $\omega$ (which is reminiscent of the fact - proved in the previous section -  that $\fm_{s,t}^\omega(m)$, $m \in \mathcal M_\lambda$, has Lebesgue measure as a `spatial' marginal for any $s$ and $t$).

Let us now define a Poisson point process $\mathcal N$ on $\T\times \mathcal{D}\times [0,\infty)$ with intensity measure $dx\, \rmPo_x(d\zeta)\otimes d\ell$. That is, we fix a Poisson point process $(x_j)_{j\in J}$ of locations at time $0$, and launch a path $\zeta^j$ from the point $x_j$. In addition, each path is labeled by a \emph{level} in $[0,\infty)$ which we shall use to indicate who reproduces during an event. Note that, by the invariance of Lebesgue measure under the dynamics of $\zeta$, the spatial distribution of the population is conserved at any time. That is:
\begin{lemma} \label{lem: PPP conserved}
For every $t\geq 0$ and every $c\in [0,\infty]$, the set $\{(\zeta_t^j,\ell_j) \, : \, (x_j,\zeta^j,\ell_j) \in \mathcal N,\, \ell_j \le c\}$ forms a Poisson point process on $\T \times [0,c]$ with intensity measure $dx \otimes \ind_{\{\ell\leq c\}}d\ell$.
\end{lemma}
Since the epochs of a reproduction event are deterministic when $\omega$ is fixed, we can in particular use Lemma~\ref{lem: PPP conserved} at the time $t_i$ of an event.

The process $\mathcal N$ encodes which individual takes part in a given reproduction event: if the time of this event is $t_i\geq 0$, then all individuals such that $\zeta^j_{t_i-}\neq \zeta^j_{t_i}$ are  affected by the event, while the others remain unaffected.
Let us now define the ancestral lineage of the individual $(t, \zeta^j_t)$ by first tracing back  $\zeta^j$ to its most recent jump before $t$. If this jump happened at time  $t_{i*}$, say, and if  $j_*$ is the index of the path that had the lowest level among all those affected at time $t_{i*}$, then the ancestral location remains at $\zeta^j_t$ for $0\le h <t-t_{i*}$ and jumps to  $\zeta^{j_*}_{t_{i*}}$ at time $h= t-t_{i*}$. Back from $t_{i*}$, we then proceed inductively.

As we did for the ancestral process $\cA$ (defined in Section \ref{section:introductory}), for any finite set of distinct indices $j_1,\ldots,j_k$ we call $\cP^t_s:=\cP^t_s(j_1,\ldots,j_k)$ the marked partition defined as
$$
\cP^t_h = \big\{(\pi^1_h,X^1_h),\ldots,(\pi^{L_h}_h,X^{L_h}_h)\big\},
$$
where each $\pi^\iota_h$ contains the labels of all the individuals in $\{j_1,\ldots,j_k\}$ at time $t$ who share a common ancestor at time $t-h$, the second component $X^\iota_h$ stands for the spatial location of this ancestor at time $t-h$, and $L_h$ denotes the number of distinct ancestors at that time. Notice that, since the Poisson point process on $\T\times \mathcal{D}\times [0,\infty)$  has an infinite intensity measure $dx\, \rmP^{\omega}_x \otimes d\ell$, it always contains infinitely many individuals even though $\T$ may have finite volume.

The following lemma is crucial for the look-down representation of the SLFV.
\begin{lemma}\label{lem: look-down}
For any finite set $\{j_1,\ldots,j_k\}$, $(\cP^t_h)_{h\in [0,t]}$ has the same law as the ancestral process $(\cA_h)_{h\in [0,t]}$ under $\rmP^{\omega,t}_{\wp_k(\zeta_t^{j_1},\ldots,\zeta_t^{j_k})}$ (recall the notation $\wp_k$ from (\ref{def pn})).
\end{lemma}
\noindent {\bf Proof.} By Lemma~\ref{lem: PPP conserved}, at any time $t\geq 0$ the set $\{(\zeta^j_t,\ell_j)\, :\,  (x_j,\zeta^j,\ell_j) \in \mathcal N\}$ forms a Poisson point process on $\T\times [0,\infty)$ with intensity measure $dx\otimes d\ell$. Furthermore, the environment $\omega$ is fixed and so the times $t_i$ at which the events take place are deterministic. For these two reasons, we first claim that
\begin{claim}\label{location anc}
During any given event, the spatial location just before the event of the affected individual with lowest level is uniformly distributed over the range of the event.
\end{claim}
Indeed, let us write $(t_i,z_i,r_i,u_i)$ for this event. By the thinning property of Poisson point processes, the set of affected individuals forms a Poisson point process on $\mathbb B(z_i,r_i)\times [0,\infty)$ with intensity measure $u_i\ind_{\mathbb B(z_i,r_i)}(x)\, dx\otimes d\ell$, and thus levels and locations are attributed independently. The result is then straightforward.

Second, let us define $\tau$ as the quantity in $[0,t]$ such that $t-\tau$ is the most recent time at which one of our $k$ paths jumps before time $t$; the time $t-\tau$ is necessarily the epoch $t_i$ of a reproduction event. Let us show that $\tau$ has the same distribution as the first time at which $\cA$ jumps under $\rmP^{\omega,t}_{\wp_k}$, where we write $\wp_k$ instead of $\wp_k(\zeta^{j_1}_t,\ldots,\zeta^{j_k}_t)$ to ease the notation. For $\nu=1,\ldots, k$, let $\alpha_\nu(t_i)$ stand for the set of indices of the events of $\omega$ occurring in the time interval $(t_i,t]$ and whose range overlaps $\zeta^{j_\nu}_t$. Let us now observe that, by independence of the paths $\zeta^j$ (which jump or remain unaffected independently of each other during an event), we can write
\begin{align}
\rmP^{\omega}[t-\tau \leq t_i]& = \rmP^{\omega}\big[\hbox{none of the $k$ paths is affected during }(t_i,t]\big] \nonumber \\
& = \prod_{\nu=1}^k \rmP^{\omega}\big[\zeta^{j_\nu}_s = \zeta^{j_\nu}_t,\ \forall\, s\in (t_i,t] \big] \nonumber \\
& = \prod_{\nu=1}^k  \prod_{a \in \alpha_\nu(t_i)}(1-u_a). \label{first jump}
\end{align}
But coming back to the definition of $\cA$ under $\rmP^{\omega,t}$, we see that (\ref{first jump}) is exactly the probability that none of the lineages starting from locations $\zeta^{j_1}_t,\ldots,\zeta^{j_k}_t$ jumps before time $t-t_i$. Hence, the first time $\tau$ at which $\cP$ jumps has the same law as the first time at which $\cA$ jumps.

Thirdly, again by the independence of $\zeta^{j_1},\ldots,\zeta^{j_k}$, all individuals present in the area of an event just before (and therefore just after) the event occurs choose independently of each other whether they take part in the event and jump, or not. Therefore, during an event affecting at least one of them, the choice of who is affected or not is decided by a set of independent Bernoulli r.v.'s, as in the evolution of $\cA$. By Claim~\ref{location anc}, the location of the individual with lowest level, onto which all affected individuals \emph{look down} to find an ancestor, is uniformly distributed over the area of the event. Furthermore, the location of an affected individual is resampled independently of its current position during the event, so that the parental location is independent of the position of its offspring. This is precisely what happens to $\cA$ during a merger.

Lastly, one can pursue the analysis by defining the most recent time $t-\tau_2$ such that the paths $\zeta^{j_1},\ldots,\zeta^{j_\nu}$ do not jump during the time interval $(t-\tau_2,t-\tau)$. The same reasoning gives us that $\tau_2$ has the same distribution as the epoch of the second jump of $\cA$ under $\rmP^{\omega,t}_{\wp_k}$, that the set of individuals who look down on the lowest level during the event occurring at time $t-\tau_2$ has the same law as the second merger of $\cA$ and that the location of the ancestor is uniformly distributed over the area of the event. Carrying on in this way and using the facts that both $\cP$ and $\cA$ are finite-rate jump processes, we can conclude that they are equal in distribution. \hfill $\Box$

\subsection{Look-down representation of the SLFV-process} \label{subs: particle approx}
As in Section~\ref{section:exist}, let us now consider a Feller process $(\cK_t)_{t\geq 0}$ with values in the compact type space $\mathbb K$. Let $B$ stand for the generator of $\mathcal{K}$. For convenience, we shall see $B$ as an operator on $C_c(\T\times \mathbb K)$ that acts only on the second coordinate. With this in mind, from now on we assume that the domain $\mathcal{D}(B)$ of $B$ is dense in $C_c(\T \times \mathbb K)$.

Using the particle system together with its genealogy from the previous subsection, for any measure $m \in \mathcal M_\lambda$, we construct a configuration of $\mathbb T \times \mathbb K$-valued paths $(\zeta^j_t, \mathcal K^j_t, \ell_j)_{j \in J}$ by choosing the $\mathcal K^j_0$ independently with distribution $\rho(x_j, \cdot)$, $j \in J$, and letting the $ \mathcal K^j$ evolve according to the look-down genealogy of the particle system. That is, mutations occur independently on each level and whenever a path $j$ looks down on a path $i$ with lower level, $\mathcal K^j$ jumps to the current value of $\mathcal K^i$ at that time.

For any $n\in \N$ and $t \ge 0$, let us define
$$
M^n_t := \frac 1n \sum_{j: \ell_j \le n} \delta_{(\zeta^j_t, \mathcal K^j_t)}.
$$
From Lemma \ref{lem: look-down} and the fact that at any time, the set of locations of the particles forms a Poisson point process with Lebesgue intensity, we obtain that for  $\P$-almost every $\omega$
$$nM_t^n \stackrel d = \mathcal N^n$$
under
 $\rmP^{\omega}$,
where the marked point configuration  $\mathcal N^n$ was defined in the proof of Theorem 1'. By an argument involving exchangeability and thinning, this extends, for any fixed time $t\ge 0$ and all $n \in \mathbb N$, to the distributional equality
$$(M_t^1, 2M_t^2, \ldots, nM_t^n) \stackrel d =  (\mathcal N^1, \mathcal N^2, \ldots ,  \mathcal N^n).$$
Hence, for any $t\ge 0$, the sequence $(M_t^n)_{n\ge 0}$ has the same distribution as the sequence  $\frac 1n \mathcal N^n$, which, as was shown in the proof of Theorem 1', converges  $\rmP^{\omega}$ a.s. for $\P$-almost every $\omega$.
Thus, for every $f\in \mathcal{D}(B)\cap C_c(\T\times \mathbb K)$, one can define the a.s. limit
\begin{equation}\label{def mi}
\la \mi_t,f\ra := \lim_{n\rightarrow \infty} \la M^n_t,f\ra = \lim_{n\rightarrow \infty}\frac{1}{n}\ \sum_{j:\ell_j\leq n}f(\zeta^j_t,\mathcal{K}^j_t).
\end{equation}
This convergence holds in fact in a pathwise manner, as shown by the following result. Recall from Lemma~\ref{lem: metric} the topology of uniform convergence over compact time intervals with which $D_{\M_\lambda}[0,\infty)$ is equipped.

\begin{thm2}
For $\P$-a.e. environment $\omega$, the sequence $(M^n)_{n\geq 1}$ converges $\rmP^{\omega}$-a.s. towards the process $(\mi_t)_{t\geq 0}$, uniformly over compact time intervals. Furthermore, $\mi$ has the same law as the quenched spatial $\Lambda$-Fleming-Viot process of Theorem~1 with initial condition $M_0=m$.
\end{thm2}

Our proof of Theorem~2 will be guided by  the proof of Theorem~1.1 in \cite{BB+2009}. The new difficulty which arises in our setting is due to the need to control the number and the spatial distribution of all the particles with levels less than $n$ lying within the compact support $S_f$ (for the `geographical' coordinate)  of some function $f\in C_c(\T\times \mathbb K)$ in a given (deterministic or random) time interval. This is achieved through the following lemma.

\begin{lemma}\label{lem: control passage}
For $\P$-almost all $\omega\in \Omega$ and every $T\geq 0$, let $\theta = \theta(\omega,T)$ be defined by
$$
\theta:=\int_{\T}\rmP^{\omega}_x\big((\zeta_t)_{0\leq t\leq T} \cap S_f \neq \emptyset\big)\,  dx.
$$
Then, $\theta$ is finite with $\P$-probability $1$.
\end{lemma}

\noindent {\bf Proof of Lemma~\ref{lem: control passage}} This property is obvious when $\T$ has finite Lebesgue measure, and so let us concentrate on the case where $\T$ is unbounded. By Fubini's theorem, we have
$$
\E[\theta]  = \int_\T \E\big[\rmP^{\omega}_x\big((\zeta_t)_{0\leq t\leq T} \cap S_f \neq \emptyset\big)\big]\, dx = \int_\T \rmP_x\big[(\zeta_t)_{0\leq t\leq T} \cap S_f \neq \emptyset\big]\, dx,
$$
where (in analogy to the notation introduced in Remark~\ref{rem: E}) $\rmP_x:=  \int \rmP^{\omega}_x(.)\P(d\omega)$ denotes the \emph{annealed} probability measure of $\zeta$. But under  $\rmP_x$, $\zeta$ is a finite-rate jump process whose instantaneous jump rates are bounded uniformly in the location by the quantity in (\ref{cond def}). Let us assume that the jump rate $\mathcal{J}$ is independent of the location and equal to this quantity (as, e.g., when $\T=\R^d$). Indeed, the existence of  boundaries only slow down the evolution of a lineage by making the balls $\bB(x,r)\subset \T$ smaller than the $d$-dimensional ball $\bB(x,r)\subset \R^d$, and so the desired property will remain satisfied even when the jump rate of $\zeta$ is inhomogeneous in space. To simplify the notation, let us also suppose that $S_f= \bB(0,a)$ for some $a>0$. Again, we do not loose generality with this assumption since there exists $a>0$ such that $S_f\subset\bB(0,a)$, and the probability of entering $S_f$ before time $T$ is bounded by that of entering $\bB(0,a)$. Since the volume of possible centres for an event of radius $r$ overlapping both $x\in \T$ and a subset of $S_f$ is bounded by the volume $V_r$ of the ball $\bB(x,r)$, we can write that
$$
\rmP_x\big[\zeta \hbox{ lies in }S_f \hbox{ after its first jump}\big] \leq \frac{1}{\mathcal{J}}\, \int_{\frac{|x|-a}{2}\vee 0}^\infty V_r\bigg(\int_0^1 u \nu_r(du)\bigg) \frac{\mathrm{Vol}(S_f)}{V_r}\, \mu(dr).
$$
(Note that a ball overlapping both $x$ and $S_f$ has radius at least $(|x|-a)/2$, then the volume of possible centres is bounded by $V_r$, the lineage jumps with probability $u$ chosen according to $\nu_r(du)$ and finally the probability that it jumps into $S_f$ is bounded by $\mathrm{Vol}(S_f)/V_r$). For $x\in \bB(0,3a)$, we shall simply bound this probability by $1$. Integrating against Lebesgue measure, turning to polar coordinates, and finally using Fubini's theorem, we obtain that
\begin{align*}
 \int_\T \rmP_x\big[\zeta &\hbox{ lies in }S_f \hbox{ after its first jump}\big]\, dx\\
& \leq \mathrm{Vol}(\bB(0,3a)) + \frac{\hbox{Vol}(S_f)C_d}{\mathcal{J}}\, \int_{3a}^\infty R^{d-1} \int_{\frac{R-a}{2}}^\infty \bigg(\int_0^1 u\nu_r(du)\bigg)\, \mu(dr)dR \\
&  \leq  \mathrm{Vol}(\bB(0,3a)) + \frac{\hbox{Vol}(S_f)C_d}{\mathcal{J}} \int_a^\infty \bigg(\int_0^1 u\nu_r(du)\bigg) \int_0^{2r+a}R^{d-1}dR\, \mu(dr) \\
&  \leq \mathrm{Vol}(\bB(0,3a)) +  C \int_a^\infty r^d \int_0^1 u\nu_r(du)\mu(dr) := \hat{C}<\infty
\end{align*}
by (\ref{cond def}), where $C_d$ is a constant depending only on the dimension $d$, and $C,\hat{C}>0$ depend on $d$ and $f$.

As concerns the subsequent jumps, recall from the beginning of Section~\ref{subs: ld} that Lebesgue measure is invariant under the evolution of $\zeta$. As a consequence, calling $\tau$ the first time at which $\zeta$ jumps and using the strong Markov property of $\zeta$, we have that
\begin{align*}
\int_\T \rmP_x \big[\zeta \hbox{ lies in }S_f \hbox{ after its 2nd jump}\big]\, dx & = \int_\T \rmE_x\big[\rmP_{\zeta_\tau}\big[\zeta \hbox{ lies in }S_f \hbox{ after its first jump}\big]\big]\, dx \\
& = \int_{\T} \rmP_x\big[\zeta \hbox{ lies in }S_f \hbox{ after its first jump}\big]\, dx \leq \hat{C}.
\end{align*}
Proceeding by recursion, we obtain that the integral over $\T$ of the probability that $\zeta$ belongs to $S_f$ after its $i$-th jump is bounded by $\hat{C}$ for any $i\geq 1$. But by our assumptions, the number of jumps made by $\zeta$ in the time interval $[0,T]$ is independent of its starting point and is a Poisson r.v. with parameter $\mathcal{J}T$. Hence, since the number of jumps of $\zeta$ is also independent of its trajectory we can write that
\begin{align*}
\E[\theta]& = \int_\T \rmP_x\big[(\zeta_t)_{0\leq t\leq T} \cap S_f \neq \emptyset\big]\, dx \\
& = \int_\T\sum_{k=0}^\infty \rmP_x\big[k \hbox{ jumps in }[0,T]\big] \rmP_x\big[\zeta\hbox{ belongs to }S_f\hbox{ initially or after one of its first }k\hbox{ jumps}\big]dx\\
& \leq  e^{-\mathcal{J}T}\sum_{k=0}^\infty \frac{(\mathcal{J}T)^k}{k!} \sum_{i=0}^k\int_\T \rmP_x\big[\zeta \hbox{ lies in }S_f\hbox{ after its }i\hbox{th jump}\big]\, dx \\
& \leq e^{-\mathcal{J}T}\sum_{k=0}^\infty \frac{(\mathcal{J}T)^k}{k!}\, \big(\hbox{Vol}(S_f)+k\hat{C}\big) = \hbox{Vol}(S_f) + \hat{C}\mathcal{J}T <\infty.
\end{align*}
We can thus conclude that $\theta<\infty$ for a.e. environment.  \hfill $\Box$

\smallskip
We can now turn to the proof of Theorem~2.

\smallskip
\noindent {\bf Proof of Theorem~2.}
Let us start by the first statement.  The key ingredient of our proof is adapted from Lemma~3.2 in \cite{BB+2009}:
\begin{lemma}\label{lem: tightness} For $\mathbb P$-a.e. environment $\omega$, for
every $T,\e>0$ and every $f\in \mathcal{D}(B)\cap C_c(\T\times \mathbb K)$, there exists a summable sequence $(\delta_n)_{n\geq 1}$ such that for every $n\in \N$,
$$
\rmP^{\omega}\bigg(\sup_{0\leq t\leq T} \big|\la M^n_t,f\ra - \la \mi_t,f\ra\big| \geq 11\e\bigg) \leq \delta_n.
$$
\end{lemma}
Indeed, let us fix $f\in \mathcal{D}(B)\cap C_c(\T\times \mathbb K)$ and for every $s\geq 0$, let $\Upsilon_s$ be defined by
\begin{equation}\label{def U}
\Upsilon_s = \sum_{(t,z,r,u)\in \omega} \ind_{\{t\leq s\}} \, u \hbox{Vol}\big(\bB(z,r)\cap S_f\big),
\end{equation}
where here again $S_f \subset \T$ stands for the compact support of $f$ (for the `geographical' coordinate). This quantity will give us an upper bound on the sum of the jumps of $\la M^n_t, f\ra$ over the time interval $[0,T]$. Note that
\begin{align}
\E[\Upsilon_T] &= \int_0^T \int_0^{\infty}\int_{S_f+\bB(0,r)}\int_0^1 u \hbox{Vol}\big(\bB(z,r)\cap S_f\big) \nu_r(du)dz\mu(dr)dt \nonumber\\
&\leq T \int_0^{\infty}\int_{S_f+\bB(0,r)}\int_0^1 u \big(C_1r^d \wedge \hbox{Vol}(S_f)\big)\nu_r(du)dz\mu(dr) \nonumber\\
& \leq T C_2 \int_0^{C_3}\int_0^1 ur^d \nu_r(du)\mu(dr) + TC_4 \int_{C_3}^{\infty}\int_0^1 u \mathrm{Vol}\big(S_f+\bB(0,r)\big) \nu_r(du)\mu(dr) \nonumber \\
& \leq T C_5 \int_0^\infty \int_0^1ur^d \nu_r(du)\mu(dr) <\infty \label{bound U}
\end{align}
for some constants $C_1,C_2,C_3,C_4, C_5>0$, where the last line uses the integrability condition (\ref{cond def}). As a consequence, with $\P$-probability $1$ the quantity $\Upsilon_T$ is finite.

Now, let us define the deterministic times $(\alpha_i)_{i\geq 1}$ by
$$
\alpha_1 := \inf\bigg\{t:\, \Upsilon_t > \frac{1}{n^2} \bigg\} \wedge \frac{1}{n^2}
$$
and
$$
\alpha_{i+1}:= \inf\bigg\{t:\, \Upsilon_t > \Upsilon_{\alpha_i}+\frac{1}{n^2} \bigg\} \wedge \bigg(\alpha_i+ \frac{1}{n^2}\bigg).
$$
This sequence of times decomposes the interval $[0,T]$ into at most $i_n:=2(\Upsilon_T+T)n^2$ subintervals, over which we shall control the fluctuations of $M^n$ and $\mi$. To this end, let us define a sequence $(\tilde{\alpha_i})_{i\geq 1}$ of stopping times (with respect to the filtration $\{\sigma(\{\mi_s,\, s\leq t\}),\, t\geq 0\}$) by
$$
\tilde{\alpha}_i := \inf\big\{t>\alpha_i:\ |\la \mi_t,f\ra - \la \mi_{\alpha_i}, f\ra| \geq 6\e\big\},
$$
where by convention $\inf \emptyset = \infty$. Let us set
$$
H_i:= |\la \mi_{\alpha_i},f\ra - \la M^n_{\alpha_i},f\ra| \vee |\la \mi_{\tilde{\alpha}_i},f\ra - \la M^n_{\tilde{\alpha}_i},f\ra|.
$$
Recall the quantity $\theta= \theta(\omega,T)$ defined in Lemma~\ref{lem: control passage}, and that the number of paths with level less than $n$ and passing through $S_f$ within the time interval $[0,T]$ is a Poisson r.v. with parameter $n\theta$. By Lemma~\ref{lem: control passage}, $\theta<\infty$ for $\P$-a.e. environments. In addition, $\theta> \mathrm{Vol}(S_f)>0$ $\P$-a.s., which guarantees that $n\theta$ is a.s. of the same order as $n$.

Let us now observe that, since $f(\zeta,k)=0$ for any $\zeta\notin S_f$, for every $t\in [0,T]$ we can write
$$
\la M^n_t,f\ra = \frac{1}{n}\, \sum_{j\in \Theta_n} f\big(\zeta^j_t,\mathcal{K}^j_t\big) = \frac{\mathrm{Card}(\Theta_n)}{n}\, \frac{1}{\mathrm{Card}(\Theta_n)}\, \sum_{j\in \Theta_n} f\big(\zeta^j_t,\mathcal{K}^j_t\big),
$$
where $\Theta_n$ is the set of all indices $j$ such that $\ell_j\leq n$ and $(\zeta^j_t)_{0\leq t\leq T}\cap S_f\neq \emptyset$. The analysis of the parameter $n \theta$ made just above shows that, if we define the event $A_n$ by
\begin{equation}\label{def ai}
A_n:=\bigg\{\hbox{Card}(\Theta_n) \notin \bigg[\frac{n\theta}{2},2 n\theta\bigg]\bigg\},
\end{equation}
then there exists $C_6>0$ (independent of $\theta$) such that for every $n\in \N$, $\rmP^\omega(A_n)\leq e^{-C_6\theta n}$. But since $\theta >0$, Lemma 3.1 in \cite{BB+2009} and the fact that
$\mathrm{Card}(\Theta_n)/n$ becomes concentrated around $\theta$ yield
\begin{equation} \label{point tight}
\rmP^{\omega}\bigg(\max_{i\leq i_n} H_i \geq \e\bigg) \leq \rmP^\omega(A_n)+\sum_{i=1}^{i_n}\rmP^{\omega}(H_i\geq \e \hbox{ and }A_n^c) \leq e^{-C_6n}+ 16(\Upsilon_T+T)n^2e^{-C_7n}
\end{equation}
for some constant $C_7>0$ depending only on $f$ and $\e$. This controls the distance between $M^n$ and $\mi$ at some discrete times. Since the definition of $\tilde{\alpha}_i$ bounds the variations of $\la \mi,f\ra$ over the small time intervals of interest, there remains to show that $\la M^n,f\ra$ does not fluctuate too much over these intervals.

For any $i\in \{1,\ldots,i_n\}$, and any index $j\in \N$, let $\tau_{i,j}$ be the first time after $\alpha_i$ at which the type and location of the individual at level $\ell_j$ is \emph{updated} during a reproduction event, with the result that the value of $f(\zeta^j,\mathcal{K}^j)$ changes (so that an individual outside $S_f$ may have been updated several times before the time $\tau_{i,j}$ at which it enters $S_f$). Splitting the evolution of the individuals of $M^n$ into its two components, mutation and reproduction, we can write that for any $t\in [\alpha_i,\alpha_{i+1}]$
\begin{align}
\la M^n_t,f\ra - \la M^n_{\alpha_i},f\ra & = \frac{1}{n}\sum_{j:\ell_j\leq n} \big\{f(\zeta^j_{(t\wedge \tau_{i,j})-},\mathcal{K}^j_{(t\wedge \tau_{i,j})-})- f(\zeta^j_{\alpha_i},\mathcal{K}^j_{\alpha_i})\big\} \nonumber\\
& \qquad + \frac{1}{n} \sum_{j:\ell_j \leq n} \big\{f(\zeta^j_t,\mathcal{K}^j_t) - f(\zeta^j_{(t\wedge \tau_{i,j})-},\mathcal{K}^j_{(t\wedge \tau_{i,j})-})\big\}. \label{fluctuations}
\end{align}
The second term in the r.h.s. of (\ref{fluctuations}) is bounded by $\frac{2\|f\|_{\infty}}{n}\mathcal{N}_n(\alpha_i,\alpha_{i+1})$, where $\mathcal{N}_n(\alpha_i,\alpha_{i+1})$ denotes the number of particles $j$ living at time $\alpha_i$, with level at most $n$ and such that $\tau_{i,j}\leq \alpha_{i+1}$. These particles are of three kinds:
\begin{itemize}
\item[$(a)$] either they belong to $S_f$ at time $\alpha_i$ and are affected by an event $(s,z,r,u)\in \omega$ before time $\alpha_{i+1}$ (each with probability $u$ during this event),
\item[$(b)$] or they do not belong to $S_f$ at time $\alpha_i$ but the first event $(t,z,r,u)$ that affects them happens before time $\alpha_{i+1}$ and brings them within $S_f$ (each with probability $u\mathrm{Vol}(S_f \cap \bB(z,r))/\mathrm{Vol}(\bB(z,r))$ during this event),
\item[$(c)$] or else they do not belong to $S_f$ at time $\alpha_i$, they jump at least twice between $\alpha_i$ and $\alpha_{i+1}$ and one of this jump (except the first one) brings them into $S_f$.
\end{itemize}
We thus need to bound the number of each kind of particles, which we denote respectively by $\mathcal{N}_{n,i}^{(a)}$, $\mathcal{N}_{n,i}^{(b)}$ and $\mathcal{N}_{n,i}^{(c)}$. To simplify the notation, let us also set $\e':= \e/(2\|f\|_{\infty})$.

First, by Lemma~\ref{lem: PPP conserved} the collection $\{(\zeta^j_{\alpha_i},\ell_j),\, j\in J\}$ forms a Poisson point process on $\T\times [0,\infty)$ with intensity measure $dx\otimes d\ell$. Hence, if all the particles initially (i.e., at time $\alpha_i$) contained in $S_f$ were frozen until an event $(s,z,r,u)$, the number of those affected by this event and with level less than $n$ would be a Poisson r.v. $\mathcal{P}_s$ with parameter $nu\mathrm{Vol}(S_f\cap \bB(z,r))$. But some of these particles may have already been affected by a previous event, so that in fact $\mathcal{P}_s$ gives a stochastic upper bound on the number of particles affected for the first time and by this event. Considering now all the events occurring between $\alpha_i$ and $\alpha_{i+1}$, we obtain that
\begin{align*}
\mathcal{N}_{n,i}^{(a)} \preceq \sum_{(s,z,r,u):\alpha_i< s\leq \alpha_{i+1}}\mathcal{P}_s  & \stackrel{(d)}{=} \mathrm{Poisson}\bigg(\sum_{(s,z,r,u):\alpha_i< s\leq \alpha_{i+1}}nu\mathrm{Vol}(S_f\cap \bB(z,r))\bigg) \\
 & = \mathrm{Poisson}\big(n (\Upsilon_{\alpha_{i+1}}-\Upsilon_{\alpha_i})\big),
\end{align*}
where the equality in distribution comes from the fact that the $\mathcal{P}_s$'s can be chosen independent since they correspond to thinnings of the set $\{(\zeta^j_{\alpha_i},\ell_j): \zeta^j_{\alpha_i}\in S_f,\, \ell_j\leq n\}$ that are independent of each others. Now, by construction $\Upsilon_{\alpha_{i+1}}-\Upsilon_{\alpha_i}\leq 1/n^2$, and so there exists a constant $C_8>0$ such that
\begin{equation}\label{particles i}
\rmP^{\omega}\big(\mathcal{N}_{n,i}^{(a)} \geq n \e'/3\big) \leq \hbox{Prob}\big(\hbox{Poisson}(1/n) \geq n\e' /3) \leq  e^{-C_8 n}.
\end{equation}

Using the same reasoning, the number of particles originally outside $S_f$ and affected for the first time during an event $(s,z,r,u)$ that brings them within $S_f$ is stochastically bounded by a Poisson r.v. with parameter
$$
nu \mathrm{Vol}(\bB(z,r)\setminus S_f)) \frac{\mathrm{Vol}(S_f \cap \bB(z,r))}{\mathrm{Vol}(\bB(z,r))} \leq nu\mathrm{Vol}(S_f\cap \bB(z,r)).
$$
Hence, here again we have that
\begin{equation}\label{particles ii}
\rmP^{\omega}\big(\mathcal{N}_{n,i}^{(b)} \geq n \e'/3\big) \leq  e^{-C_8 n}.
\end{equation}

Finally, $\mathcal{N}_{n,i}^{(c)}$ is stochastically bounded by a Poisson r.v. with parameter
\begin{equation}\label{double jumps}
n\int_{\T}\rmP^{\omega}_x\big(\hbox{at least 2 jumps in }[\alpha_i,\alpha_i+n^{-2}] \hbox{ and passes through }S_f \hbox{ in this interval}\big)\, dx.
\end{equation}
Observe that $\alpha_i$ is a stopping time with respect to the Poisson point process of reproduction events. Hence, using the Markov inequality and then Fubini's theorem, we can write that
\begin{align*}
\P\bigg[ &\int_{\T}\rmP^{\omega}_x\big(\hbox{at least 2 jumps in }[\alpha_i,\alpha_i+n^{-2}] \hbox{ and passes through }S_f \hbox{ in this interval}\big)\, dx\geq \e'/6\bigg] \\
&\leq \frac{6}{\e'}\, \int_{\T}\rmP_x\big(\hbox{at least 2 jumps in }[\alpha_i,\alpha_i+n^{-2}] \hbox{ and passes through }S_f \hbox{ in this interval}\big)\, dx.
\end{align*}
Proceeding as in the proof of Lemma~\ref{lem: control passage}, we obtain that the above integral is bounded by a constant times $n^{-4}$. Consequently,
$$
\P\bigg[\exists i\leq i_n : \int_{\T}\rmP^{\omega}_x\big(\geq \hbox{2 jumps in }[\alpha_i,\alpha_i+n^{-2}] \hbox{ and passes through }S_f\\
  \hbox{ in this interval}\big)\, dx\geq \e'/6\bigg]  \leq \frac{C_9}{n^2}.
$$
(Recall that $\E[i_n]$ is proportional to $n^2$.) By the Borel-Cantelli lemma, we can conclude that for $\P$-a.e. environments, there exists $n_0(\omega)$ such that for every $n\geq n_0$,
$$
\int_{\T}\rmP^{\omega}_x\big(\geq \hbox{2 jumps in }[\alpha_i,\alpha_i+n^{-2}] \hbox{ and passes through }S_f \hbox{ in this interval}\big)\, dx< \e'/6, \qquad \forall i\leq i_n.
$$
Assuming that our fixed environment belongs to this set of full $\P$-probability and recalling (\ref{double jumps}), we obtain that for every $n\geq n_0(\omega)$
\begin{equation}\label{particles iii}
\rmP^\omega\big(\exists i\leq i_n\, :\, \mathcal{N}_{n,i}^{(c)}\geq n\e'/3\big)\leq C_{10} n^2 e^{-C_{11}n}.
\end{equation}
Combining (\ref{particles i}), (\ref{particles ii}), (\ref{particles iii}) and the facts that $2\|f\|_{\infty}\e'=\e$ and
$$
\mathcal{N}_n(\alpha_i,\alpha_{i+1}) = \mathcal{N}_{n,i}^{(a)} + \mathcal{N}_{n,i}^{(b)} + \mathcal{N}_{n,i}^{(c)},
$$
we obtain that
\begin{equation}
\rmP^{\omega}\bigg(\max_{i\leq i_n}\sup_{\alpha_i<t\leq \alpha_{i+1}}\bigg|\frac{1}{n} \sum_{j:\ell_j \leq n} \big\{f(\zeta^j_t,\mathcal{K}^j_t) - f(\zeta^j_{(t\wedge \tau_{i,j})-},\mathcal{K}^j_{(t\wedge \tau_{i,j})-})\big\}\bigg|\geq \e\bigg) \leq C_{12} n^2 e^{-C_{13}n} \label{control 1}
\end{equation}
for some constants $C_{12},C_{13}>0$.

Let us now consider the first term in the r.h.s. of (\ref{fluctuations}), which corresponds to the evolution due to mutation. Recall that $B$ denotes the generator of the mutation process, and that mutations occur independently along distinct lines. Recall also that we write $Bf$ for the action of $B$ on the second argument of $f$. We have
\begin{align*}
\frac{1}{n}\sum_{j:\ell_j\leq n} &\big\{f\big(\zeta^j_{(t\wedge \tau_{i,j})-},\mathcal{K}^j_{(t\wedge \tau_{i,j})-}\big)- f\big(\zeta^j_{\alpha_i},\mathcal{K}^j_{\alpha_i}\big)\big\} \\
& = \frac 1n \sum_{j:\ell_j\leq n} \bigg\{f\big(\zeta^j_{(t\wedge \tau_{i,j})-},\mathcal{K}^j_{(t\wedge \tau_{i,j})-}\big)- f\big(\zeta^j_{\alpha_i},\mathcal{K}^j_{\alpha_i}\big) - \int_{\alpha_i}^{t\wedge \tau_{i,j}} Bf\big(\zeta^j_s,\mathcal{K}^j_s\big)\, ds\bigg\} \\
& \qquad \qquad + \frac 1n \sum_{j:\ell_j\leq n} \int_{\alpha_i}^{t\wedge \tau_{i,j}} Bf\big(\zeta^j_s,\mathcal{K}^j_s\big)\, ds \\
& := \Delta_i(t)+  \frac 1n \sum_{j:\ell_j\leq n} \int_{\alpha_i}^{t\wedge \tau_{i,j}} Bf\big(\zeta^j_s,\mathcal{K}^j_s\big)\, ds,
\end{align*}
where $(\Delta_i(t))_{t\in [\alpha_i,\alpha_{i+1}]}$ is a zero-mean martingale (as the sum of finitely many zero-mean martingales). Recall the event $A_n$ defined in (\ref{def ai}). For any $n$ large enough so that $2\theta n^{-2}\|Bf\|_{\infty} \leq \e$, we have
\begin{align*}
\rmP^{\omega}\bigg(&\max_{i\leq i_n}\sup_{\alpha_i<t\leq \alpha_{i+1}}\bigg|\frac{1}{n}\sum_{j:\ell_j\leq n} \big\{f\big(\zeta^j_{(t\wedge \tau_{i,j})-},\mathcal{K}^j_{(t\wedge \tau_{i,j})-}\big)- f\big(\zeta^j_{\alpha_i},\mathcal{K}^j_{\alpha_i}\big)\big\}\bigg|\geq 2\e\bigg) \\
& \leq \rmP^\omega(A_n) +\sum_{i\leq i_n} \rmP^{\omega}\bigg(\sup_{\alpha_i<t\leq \alpha_{i+1}}|\Delta_i(t)|+ \frac{2n\theta}{n}\, n^{-2}\|Bf\|_{\infty}\geq 2\e\, ;\, A_n^c\bigg) \\
& \leq \rmP^\omega(A_n) + \sum_{i\leq i_n} \rmP^{\omega}\bigg(\sup_{\alpha_i<t\leq \alpha_{i+1}}|\Delta_i(t)|\geq \e\, ;\, A_n^c\bigg).
\end{align*}
Hence, using Lemma~3.1$(b)$ in \cite{BB+2009} together with the same argument as in (\ref{point tight}) to control the number of individuals in $S_f$ and with label at most $n$ at any given time, we obtain that
\begin{align*}
\rmP^{\omega}\bigg(\max_{i\leq i_n}\sup_{\alpha_i<t\leq \alpha_{i+1}}\bigg|\frac{1}{n}\sum_{j:\ell_j\leq n} \big\{f\big(\zeta^j_{(t\wedge \tau_{i,j})-},\mathcal{K}^j_{(t\wedge \tau_{i,j})-}\big)- f\big(\zeta^j_{\alpha_i},\mathcal{K}^j_{\alpha_i}\big)\big\}\bigg|\geq 2\e\bigg) \leq \big(2(\Upsilon_T+T)n^2+1) e^{-C_{14}n}
\end{align*}
for some $C_{14}>0$. Together with (\ref{fluctuations}) and (\ref{control 1}), this gives us that
\begin{equation}\label{control 2}
\rmP^{\omega}\bigg(\max_{i\leq i_n}\sup_{\alpha_i<t\leq \alpha_{i+1}}\big|\la M^n_t,f\ra - \la M^n_{\alpha_i},f\ra\big|\geq 3\e\bigg) \leq C_{15}n^2 e^{-C_{16}n}
\end{equation}
for some $C_{15},C_{16}>0$.

Now, it is easy to see that if $\max_{i\leq i_n} H_i < \e$ and $\max_{i\leq i_n}\sup_{\alpha_i<t\leq \alpha_{i+1}}\big|\la M^n_t,f\ra - \la M^n_{\alpha_i},f\ra\big|< 3\e$, necessarily $\tilde{\alpha}_i > \alpha_{i+1}$ for every $i\leq i_n$, which implies that
$$
\max_{i\leq i_n} \sup_{\alpha_i\leq t\leq \alpha_{i+1}} \big|\la \mi_t,f \ra - \la \mi_{\alpha_i},f\ra\big| < 6\e
$$
on this event. Using (\ref{point tight}) and (\ref{control 2}), we finally obtain that
$$
\rmP^{\omega}\bigg(\sup_{0\leq t\leq T} \big|\la M^n_t,f\ra - \la \mi_t,f\ra\big| \geq 11\e\bigg) \leq C_{17}n^2e^{-C_{18}n} =: \delta_n,
$$
and Lemma~\ref{lem: tightness} is proved. Then, the Borel-Cantelli lemma gives us the a.s. convergence of $\la M^n,f\ra$ towards $\la \mi,f\ra$ as $n\rightarrow \infty$, uniformly over compact time intervals. Since by assumption the set $\mathcal{D}(B)\cap C_c(\T\times \mathbb K)$ is separable and is dense in $C_c(\T\times \mathbb K)$, the a.s. convergence of $M^n$ towards $\mi$ follows from Lemma~\ref{lem: metric}, as well as the fact that $\mi$ has c\`adl\`ag paths with probability $1$.

\smallskip
Finally, let us now prove that $\mi$ has the same law as the \emph{quenched} spatial $\Lambda$-Fleming-Viot process of Theorem~1. As we have already mentioned at the beginning of Section~\ref{subs: particle approx}, with the notation of the proof of Theorem~1' we have
$$
M^n_t \stackrel{(d)}{=} \frac{1}{n}\,\mathcal{N}^n.
$$
Letting $n$ tend to infinity and using the a.s. convergence of $n^{-1}\mathcal{N}^n$ towards the random measure $\fm_{0,t}$, we can conclude that
$\mi_t $ has law $Q_{0,t}^\omega (m,.)$.

Let us now fix $0<s<t$ and condition on $\mi_s=m\in \M_{\lambda}$. By construction, the types of all the individuals alive at time $t$ are determined by propagating the types of the individuals living at time $s$ along the genealogical trees created by the look-down dynamics. But Lemma~\ref{lem: look-down} and the Markov property of the genealogical processes ensure that the ancestries between times $s$ and $t$ of any finite sample of individuals living at time $t$ have the same law as the trees $(\mathcal{A}_h)_{0\leq h\leq t-s}$ under the probability measure $\rmP^{\omega,t}$ (recall that $\mathcal{A}$ is defined in the paragraph around (\ref{notation A}), and $\rmP^{\omega,t}$ is introduced just after Lemma~\ref{lem:topo}). Consequently, the proof of Theorem~1' shows that $\mi_t$ has law $Q^{\omega}_{s,t}(m,dm')$. We can therefore conclude that $\mi$ is a Markov process with transition semigroup $(Q_{s,t}^\omega)_{0\leq s\leq t}$, and so Theorem~1 guarantees that $\mi$ has the same law as the quenched SLFV $(M_t)_{t\geq 0}$. The proof of Theorem~2 is now complete. \hfill $\Box$

\medskip
Let us end this section with the proof of Lemma~\ref{lem: fv}.

\smallskip
\noindent \textbf{Proof of Lemma~\ref{lem: fv}.} Again, let us consider $f\in C_c(\T\times \mathbb K)$ and show that $(\la \mi_t,f\ra)_{t\geq 0}$ has paths of finite variation with $\rmP^{\omega}$-probability $1$. By taking a countable basis $f_1,f_2,\ldots$ of $C_c(\T\times \mathbb K)$, we shall then conclude that this property holds for all $f_i$'s simultaneously with probability $1$, and so by Lemma~\ref{lem: metric} that a.s. $\mi$ has paths of finite variation in $D_{\M_\lambda}[0,\infty)$.

Let us thus fix $T > 0$. Since there is no mutation, $\mi$ evolves only at the times of a reproduction event. By construction, during the event $(t,z,r,u)$, at each site $x$ within $\bB(z,r)$ a fraction $u$ of the population is replaced by individuals that are all of some type $\kappa$. That is (in the notation of (\ref{corr})),  $\rho^{\infty}_t(x) = (1-u)\rho^{\infty}_{t-}(x) + u\delta_\kappa$. Together with the fact that the spatial marginal of each $\mi_t$ is Lebesgue measure, for any finite number of times $t_0,\ldots,t_l$ such that $0=t_0 < t_1 < \ldots < t_l = T$, we have
\begin{align*}
\sum_{j=0}^{l-1} \big|\la \mi_{t_{j+1}},f\ra - \la \mi_{t_j},f\ra \big|& \leq \sum_{(t,z,r,u)\in \omega:0< t\leq T} \big|\la \mi_t,f\ra - \la \mi_{t-},f\ra \big| \\
&\leq 2\|f\|_{\infty}\sum_{(t,z,r,u)\in \omega:0< t\leq T} u\, \hbox{Vol}(S_f\cap \bB(z,r)),
\end{align*}
where here again $S_f$ stands for the compact support of $f$ with respect to the spatial coordinate. But by (\ref{def U}) and (\ref{bound U}), the sum in the r.h.s. above is finite for $\P$-a.e. $\omega$, and is independent of the subdivision $(t_0,\ldots,t_l)$ chosen. Since the total variation of $(\la \mi_t,f\ra)_{0\leq t\leq T}$ is given by
$$
\sup_{l\in \N} \sup_{(t_0,\ldots,t_l)} \sum_{j=0}^{l-1} \big|\la \mi_{t_{j+1}},f\ra - \la \mi_{t_j},f\ra \big|,
$$
the desired result follows. \hfill $\Box$

\subsection{Generalizations}\label{subs: generalizations}

The construction carried out in the last paragraph is  robust in the sense that it can accommodate different variants of the spatial $\Lambda$-Fleming-Viot process such as that introduced in \cite{BEK2010}. In the Gaussian model of \cite{BEK2010}, the environment is given by a Poisson point process of epochs and centres with intensity measure $c\, dt\otimes dx$, where $c>0$. During an event centered on $z\in \T$, the fraction $u(z,x)$ of the population killed at site $x$ is given by
\begin{align}\label{gauss}
u(z,x):= u_0\, \exp\big(-|z-x|^2/(2\theta^2)\big),
\end{align}
where $u_0\in (0,1]$ is the maximal killing intensity and $\theta^2>0$ is a fixed parameter. Then, the location of the parent is sampled according to the kernel
\begin{align}\label{gauss2}
v(z,y):= \frac{1}{(2\pi (\alpha \theta)^2)^{d/2}}\, \exp\big( -|z-y|^2/(2\alpha^2\theta^2)\big),
\end{align}
where $\alpha >0$ (note that if $\alpha>1$, the parent is chosen from a wider area than that over which the impact of the event is non-negligible). A type is then drawn from the local type distribution at this site. Finally, the population removed during the first step is replaced by offspring of the chosen parent so that the local mass of individuals at every site remains constant. This process can be phrased as an $\M_{\lambda}$-valued evolution, and the corresponding genealogical process can be described explicitly.

In order to adapt the look-down construction of Section~\ref{subs: ld}, let us emphasize its two main ingredients:
\begin{description}
\item[An ancestry given by a finite-rate jumps process:] A lineage at location $x$ is affected by an event whenever it belongs to the fraction $u(z_i,x)$ replaced during the event. Hence, with \eqref{gauss}, the rate at which a lineage jumps is equal to
  $$
  c \int_{\T} u(z,x)\, dz = c\, u_0 \int_{\T} e^{-|x-z|^2/(2\theta^2)}\, dz <\infty.
  $$
This fact enables us to introduce a well-defined procedure to make individual types evolve in time.
\item[Reversibility of Lebesgue measure:] With the rule \eqref{gauss2}, the transition kernel $K_i$ corresponding to the $i$-th event, centered on $z_i$, is given by
$$
K_i(x,dy)= \big(1-u(z_i,x))\,\delta_x + u(z_i,x)\,v(z_i,y)\, dy.
$$
Since conditionally on $z_i$ the new location of $\zeta$ is independent of its previous position, it is not difficult to see that
$$
dx\, K_i(x,dy)=dy\, K_i(y,dx).
$$
Hence, Lebesgue measure is again reversible for the evolution of $\zeta$.
\end{description}
Next, when $\alpha=1$ we can again choose the parent to be the affected individual with lowest level. When $\alpha\neq 1$, the kernel with which the parent is drawn differs from that with which individuals are affected, and so we have to change the way the parent is sampled. This time, we thin the Poisson point process of individuals at the time of the event by keeping each point $(\zeta^j_{t_i-},\ell_j)$ with probability $v(z_i,\zeta^j_{t_i-})/v(0,0)$. We then take the individual of the thinned point process with lowest level as the parent, and resample its location using the kernel $v(z_i,y)\,dy$. In this way, some of the affected individuals may now have a lower level than the parent. However, the resampling of the location of the parent ensures that the position of the ancestor just after a given merger of $\cP$ has density $v(z_i,\cdot\,)$, which is the essential point in the equality in distribution of $\cP$ and $\cA$ (cf. Lemma~\ref{lem: look-down}).

Using the two properties emphasized above and the ideas of Sections~\ref{subs: ld} and \ref{subs: particle approx}, one can then construct a particle representation for the stochastic flow obtained in the Gaussian model.
\\\\
Another possible generalization is to allow multiple parents. Extending the ball model from Section~\ref{subs: ld}, we may e.g. fix a distribution $\gamma$ on $\N$ with compact support and associate an independent realization $N_i$ of $\gamma$ to each event of $\omega$. During the $i$-th event, $N_i$ parents are picked independently and uniformly over the area of the event. If the corresponding types are $\kappa_1,\ldots,\kappa_{N_i}$ (not necessarily distinct), the new value of the SLFV is given by
$$
M_{t_i}:= \ind_{\mathbb B(z_i,r_i)^c}\, M_{t_i-} + \ind_{\mathbb B(z_i,r_i)}\bigg\{(1-u_i)M_{t_i-} + \frac{u_i}{N_i}\sum_{n=1}^{N_i}\delta_{\kappa_n}\bigg\}.
$$
In words, within the area of the event we keep a fraction $1-u_i$ of the population as it was just before the event, and replace a fraction $u_i$ by offspring of the $N_i$ reproducing types in equal proportions. Thus, the associated genealogies can have multiple and simultaneous mergers whenever $\gamma$ puts some mass onto $\{2,3,\,\ldots\}$.

In the corresponding look-down construction, the only difference with that of Section~\ref{subs: particle approx} is that we use the $N_i$ affected individuals with lowest levels as the parents during event $i$, and decide that an affected individual looks down onto one of these $N_i$ levels with equal probability, independently of each other. Since the two conditions given above are fulfilled (with the kernel $K_i$ given by (\ref{kernel})), we can construct a look-down coupling between the \emph{quenched} SLFV and its genealogies.

\section{Coming down from infinity}\label{section:CDI}
Recall that a coalescent is said to \emph{come down from infinity} if, starting from countably many lineages, there exists a time in the past at which the number of ancestors is finite. For non-spatial exchangeable coalescents, it is known that whenever the quantity corresponding to the impact $u$ here is always less than $1$, then either the coalescent comes down from infinity instantaneously with probability $1$, or the number of ancestors remains infinite for all times a.s. See e.g. Proposition~23 in \cite{PIT1999} for a statement of this result for coalescents with multiple mergers, and Lemma~31 in \cite{SCH2000} for the more general case of coalescents with simultaneous and multiple mergers. Furthermore, a precise criterion for CDI is obtained in \cite{SCH2000b} for $\Lambda$-coalescents, and some conditions are given in Section~5.5 of \cite{SCH2000} as concerns the $\Xi$-coalescents.

In the context of spatial coalescents, only a few results exist concerning the question of coming down from infinity. Indeed, the geographical movement of lineages may separate them before they find a common ancestor, or may bring back together some lineages which were too far away from each other to coalesce. Hence, understanding the form of the resulting genealogical process requires a fine analysis of the interplay between these two mechanisms. In \cite{LS2006}, Limic and Sturm consider a population spread over the vertices of a finite graph. They assume that the lineages migrate independently of each other between the sites and can coalesce only when they belong to the same subpopulation. They show that the corresponding spatial coalescent comes down from infinity if and only if its non-spatial counterpart does. The case of discrete but infinite graphs is then explored in \cite{ABL2010}, where it is shown that because the timescales of migration and coalescence are precisely the same, an arbitrarily large number of lineages can escape from the others without coalescing. As a consequence, the spatial coalescents they consider never come down from infinity. Finally, a version of the SLFV in which the geographical space is a self-similar hierarchical structure (such as a Cantor set, or an infinite $m$-ary tree) is introduced in \cite{FRE2011}. The forwards-in-time evolution is formulated as a stochastic flow in the same spirit as Bertoin and Le Gall's construction of the $\Lambda$-Fleming-Viot superprocess in \cite{BLG2003}. Freeman then details the five different possible behaviours of the corresponding ancestral process, making explicit the importance of the geographical structure in the form of the genealogies.

Coming back to the SLFV, imagine we sample countably many individuals in the ball $B:=\bB(0,1)\subset \T$ and trace back their genealogical relations. We want to show that, under the condition
\begin{equation}\label{cond CDI}
\int_0^{\infty} \nu_r(\{1\})\,\mu(dr) =0,
\end{equation}
in a.e. environment the number of distinct ancestors remains infinite for all times in the past.

\begin{rema}
The condition (\ref{cond CDI}) simply says that the impact parameter is always less than $1$. If it does not hold, an event overlapping all the lineages but for a finite number (say, $n$) of them may have impact $1$, in which case the number of ancestors would come down to $n+1<\infty$.
\end{rema}

In fact, as one might expect from the comparison of the finiteness condition \eqref{cond def} with the condition
$$
\hbox{Rate at which a single lineage is affected}=\int_0^1 \frac{\Lambda(du)}{u} <\infty
$$
for the presence of dust in a nonspatial $\Lambda$-coalescent (see Theorem~8 in \cite{PIT1999}), we shall show that at time $t>0$ in the past infinitely many lineages have not yet been affected by an event. The difficulty here comes from the spatial correlations between the rates at which close-by lineages are affected by reproduction events. Despite this effect of space, the proof of the following proposition relies on the Poissonian structure of the events. Recall the condition (\ref{cond def}) of existence of the spatial $\Lambda$-Fleming-Viot process on $\T$ with parameters $\mu$ and $\{\nu_r,\, r>0\}$.

\begin{prop}\label{prop:CDI} Assume that (\ref{cond CDI}) holds. Let $\cN$ be a Poisson point process on $B\times [0,\infty)$ with intensity measure $dx\big|_B\otimes d\ell$, and let us use $\{(\{i\},\, x_i):\, (x_i,\ell_i)\in \cN\}$ as the initial value of the spatial $\Lambda$-coalescent $\cA$ defined in Section~\ref{subs:previous def}. Then for almost every environment $\omega$, at any time $t>0$ the set of the ancestral locations in $\cA_t$ contains a Poisson point process on $B$ with infinite intensity. In particular, the spatial $\Lambda$-coalescent $\cA$ never comes down from infinity.
\end{prop}

\noindent{\bf Proof of Proposition~\ref{prop:CDI}.} For every $x\in B$, let us write $\alpha(x,t)$ for the set of events of $\omega$ between (backward) times $0$ and $t$ overlapping $x$. We have
\begin{equation}\label{rate x}
\P\big(\rmP^{\omega,0}_x (\xi \hbox{ does not jump until time }t)>0\big) = \P\bigg(\prod_{i\in \alpha(x,t)} (1-u_i)>0\bigg)=1.
\end{equation}
Indeed, using the exponential formula for the Poisson point process of events, we can write for every $\theta>0$
\begin{equation}\label{laplace}
\E\Big[e^{\theta \sum_{i\in \alpha(x,t)}\log(1-u_i)}\Big] = \exp\left\{t\int_0^{\infty}\int_0^1 V_{r,x}\big[(1-u)^\theta -1\big]\nu_r(du) \mu(dr)\right\},
\end{equation}
where $V_{r,x}$ stands for the volume of the ball $B(x,r)$ in $\T$. But for every $u\in [0,1)$, $1-(1-u)^\theta$ decreases to $0$ as $\theta$ decreases to $0$. Furthermore, using this property with $\theta=1$ we obtain that for every $\theta\in (0,1)$,
$$
0\leq \int_0^{\infty}\int_0^1V_{r,x}\big[1-(1-u)^\theta\big] \nu_r(du)\mu(dr) \leq \int_0^{\infty}\int_0^1 V_{r,x}u\, \nu_r(du) \mu(dr) <\infty
$$
by the condition (\ref{cond def}) imposed on $\mu$ and $\{\nu_r,\, r>0\}$. We can therefore use dominated convergence and the fact that the set $\{u=1\}$ is never charged to conclude that
$$
\lim_{\theta\rightarrow 0}\int_0^{\infty}\int_0^1 V_{r,x}\big[(1-u)^\theta -1\big] \nu_r(du)\mu(dr)= 0.
$$
Coming back to (\ref{laplace}) and letting $\theta$ tend to $0$, we obtain that $\P\big[\sum_{i\in \alpha(x,t)}\log(1-u_i)=-\infty\big]=0$ and thus (\ref{rate x}) is proved.

Using now Fubini's theorem, we can write
$$
1=\frac{1}{\hbox{Vol}(B)}\int_B \P\bigg[\prod_{i\in \alpha(x,t)}(1-u_i)>0\bigg]\, dx = \E\bigg[\frac{1}{\hbox{Vol}(B)}\int_B \ind_{\big\{\prod_{i\in \alpha(x,t)} (1-u_i)>0\big\}}\, dx\bigg],
$$
and since the quantity within the expectation on the right-hand side belongs to $[0,1]$, it therefore needs to be $1$ $\P$-a.s.

In particular, for a.e. $\omega$, we can find $\e(\omega)\in (0,1)$ such that
$$
\int_B \ind_{\big\{\prod_{i\in \alpha(x,t)} (1-u_i)>\e(\omega)\big\}}\, dx >0.
$$
Let us call $\mathcal{S}_{\e,\omega}$ the set $\big\{x\in B:\, \prod_{i\in \alpha(x,t)} (1-u_i)>\e(\omega)\big\}$. What we have just shown is that for almost every $\omega$, $\mathcal{S}_{\e,\omega}$ has positive Lebesgue measure. Hence it contains infinitely many points of $\cN$, and each of them has probability at least $\e(\omega)$ not to have been affected yet, independently of each other. Together with initial Poissonian structure of $\cN$, this completes the proof of Proposition~\ref{prop:CDI}.
\hfill $\Box$

\paragraph{Acknowledgement} We thank Peter Pfaffelhuber for a stimulating remark concerning the representation of spatial Lambda Fleming-Viot processes as measure-valued processes, and Vlada Limic for her valuable contribution to the onset of this project. We also thank Alison Etheridge and Tom Kurtz for several illuminating discussions on the SLFV and general look-down constructions. Finally, we thank the referee and the Associate Editor for their careful reading and their comments which lead to an improvement of the presentation of the paper.

\end{document}